\documentclass[onefignum,onetabnum]{siamart220329}
\allowdisplaybreaks

\usepackage{lipsum}
\usepackage[export]{adjustbox}
\usepackage{amsfonts}
\usepackage{graphicx}
\usepackage{epstopdf}
\usepackage{algorithmic}
\usepackage{cancel}
\usepackage{parskip}

\usepackage{mathtools}
\usepackage{stackengine}
\usepackage{fixmath}
\usepackage{mathrsfs}
\usepackage{dsfont}
\usepackage{amsfonts}
\usepackage{xcolor}
\usepackage{hyperref}
\usepackage[capitalize]{cleveref}
\usepackage{cases}
\usepackage{graphicx} 
\DeclareMathAlphabet{\mathbold}{OML}{cmm}{b}{it}
\mathtoolsset{centercolon}  
\usepackage{amsmath, amssymb,amsbsy}
\usepackage{bm}

\renewcommand{\phi}{\varphi}
\renewcommand{\epsilon}{\varepsilon}

\newcommand{\cnst}[1]{\mathrm{#1}}

\newcommand{\zerovct}{\vct{0}} 
\newcommand{\Id}{\mathbf{I}} 

\newcommand{\R}{\mathbb{R}}

\newcommand{\vct}[1]{\mathbold{#1}}

\newcommand{\tp}{{T}}

\newcommand{\trace}{\operatorname{tr}}

\DeclarePairedDelimiterX{\infdivx}[2]{(}{)}{%
  #1\;\delimsize\|\;#2%
}

\newcommand{\divergence}{\operatorname{div}}
\newcommand\D{\mathop{}\cnst{d}}

\newcommand{\jump}[1]{\left[\!\left[#1\right]\!\right]}
\newcommand{\avg}[1]{\left\{\!\!\left\{#1\right\}\!\!\right\}}

\newcommand{\defeq}{\vcentcolon=}

\newcommand{\sens}[1]{\hat{#1}}
\newcommand{\adj}[1]{\mathring{#1}}

\definecolor{claudeviolet}{RGB}{109, 74, 155}

\usepackage{autonum}

\usepackage{tikz}
\usepackage{pgfplots}
\usetikzlibrary{shapes.arrows, patterns, calc}
\usetikzlibrary{external}
\usepackage{tikz-3dplot}
\ifpdf
  \DeclareGraphicsExtensions{.eps,.pdf,.png,.jpg}
\else
  \DeclareGraphicsExtensions{.eps}
\fi

\usepgfplotslibrary{groupplots}
\usetikzlibrary{arrows}
\usetikzlibrary{shapes}
\usetikzlibrary{decorations.text}
\usetikzlibrary{quantikz}
\usepackage{siunitx}
\usetikzlibrary{arrows.meta}
\pgfplotsset{ compat=1.18,
    standard/.style={
    scale only axis,
    width=0.5\textwidth,
    enlarge x limits=0.05,
    enlarge y limits=0.05,
    max space between ticks=40,
    every axis/.append style={font=\small},
	every legend/.append style={font=\small},
	every node/.append style={font=\small},
	ylabel style={yshift=-1.5em},
	unbounded coords=jump,
	}
}

\definecolor{steelblue}{HTML}{A1BDC7}
\definecolor{orange}{HTML}{D98C21}
\definecolor{silver}{HTML}{B0ABA8}
\definecolor{rust}{HTML}{B8420F}
\definecolor{seagreen}{HTML}{2E6B69}
\definecolor{joshua}{HTML}{FBDC7F}
\definecolor{darksky}{HTML}{154c79}

\colorlet{lightsilver}{silver!30!white}
\colorlet{darkorange}{orange!85!black}
\colorlet{darksilver}{silver!85!black}
\colorlet{darksteelblue}{steelblue!85!black}
\colorlet{darkrust}{rust!85!black}
\colorlet{darkseagreen}{seagreen!85!black}

\hypersetup{colorlinks=true,linkcolor=darkrust,citecolor=darkseagreen,urlcolor=darksilver}

\ifpdf
  \DeclareGraphicsExtensions{.eps,.pdf,.png,.jpg}
\else
  \DeclareGraphicsExtensions{.eps}
\fi

\newsiamremark{remark}{Remark}
\newsiamremark{hypothesis}{Hypothesis}
\crefname{hypothesis}{Hypothesis}{Hypotheses}
\newsiamthm{claim}{Claim}

\headers{Information geometric regularization for computing sensitivities}{Florian Sch{\"a}fer}

\title{Information geometric regularization \\ for computing sensitivities \\of flows with shocks} 

\author{Florian Sch{\"a}fer\thanks{Courant Institute of Mathematical Sciences, New York University
  (\email{florian.schaefer@nyu.edu}).}}

\usepackage{amsopn}

\ifpdf
\hypersetup{
  pdftitle={Information geometric regularization for sensitivities of flows with shocks},
  pdfauthor={Florian Sch{\"a}fer}
}

\begin{document}

\maketitle

\begin{abstract}
     Computing (adjoint) sensitivities of flows with shocks is a longstanding problem in computational fluid dynamics. 
     The spurious sensitivities due to shock sensors and limiters frequently force practitioners to accept the errors incurred by ``freezing'' limiters and shock sensors. 
     The recently proposed information geometric regularization (IGR) is an inviscid, PDE-based regularization of the compressible Euler equations that replaces shocks with smooth profiles without damping fine-scale structures. 
     This work derives the forward and adjoint sensitivities for the IGR system with periodic boundary conditions and demonstrates their convergence, under grid refinement, to the sensitivities obtained by finite differences or automatic differentiation through the forward solve.
\end{abstract}



\section{Introduction}
\subsection{The Euler equations of gas dynamics \nopunct}
describe the dynamics of inviscid fluids and gases by the system of conservation laws 
\begin{equation}
\label{eqn:euler}
\begin{cases}
     \partial_t \vct{\mu} + \divergence \left(\vct{\mu} \vct{u}^{\tp} + P \Id \right) &= \vct{f}_{\vct{\mu}}, \\
     \partial_t \rho + \divergence\left( \rho \vct{u}^{\tp} \right) &= f_\rho, \\
     \partial_t E + \divergence ((E + P)\vct{u}^{\tp}) &= \vct{u}^{\tp} \vct{f}_{\vct{\mu}} + f_E.
\end{cases}
\end{equation}
The system governs the evolution of momentum density $\vct{\mu} = \rho \vct{u}$, mass density $\rho$, and total energy density $E$.
Here, $\vct{u} = \vct{\mu}/\rho$ is the velocity, $P = P(\rho, e)$ is the mechanical pressure depending on density and specific internal energy $e = E/\rho - \vct{u}^{\tp} \vct{u}/2$.
The forcing terms $\vct{f}_{\vct{\mu}}, f_\rho, f_E$ represent external body forces or source terms.

\subsection{Shock formation in Euler's equations}
High-Mach-number gas dynamics form \emph{shock waves}, which manifest themselves as discontinuities of the velocity, pressure, and energy of the gas on macroscopic scales.
The viscosity of physical gases, although small, ensures smooth transitions on microscopic scales.
Solutions of the Euler equations describe gases with negligible viscosity and form discontinuities in finite time.
After shock formation, solutions are defined weakly, using a suitable entropy criterion or (equivalently) by a vanishing viscosity limit \cite{leveque1992numerical}.
Shock singularities cause Gibbs--Runge oscillations in higher-order numerical methods.
In geometrically simple cases, this can be avoided by \emph{shock tracking}, splitting the computational domain across the shock.
But in most cases involving dynamically evolving shocks, this is impractically complicated.
Instead, so-called \emph{shock-capturing} methods adaptively add dissipation to regularize the simulation near the shock.
Localized diffusive regularizations modify the Euler equations with a nonlinear viscous term \cite{vonneumann1950method,puppo2004numerical,guermond2011entropy,cook2005hyperviscosity,persson2006sub,mani2009suitability,kawai2010assessment,chan2025artificial,jameson2017origins}.
Limiters instead locally reduce the order of the numerical method, regularizing the shock by numerical diffusion \cite{van1979towards,harten_tvd_1983,harten1997uniformly,shu2006essentially,jiang1996efficient,guermond2018second}.

\subsection{Sensitivities of Euler solutions, problems due to shocks}
Design optimization \cite{jameson1988aerodynamic,jameson1994control,giles2000introduction}, uncertainty quantification, and scientific machine learning \cite{bryngelson2024fast,liu2024adjoint} in gas dynamics not only require solutions to Euler's equations but also their sensitivities (derivatives) with respect to problem parameters.
These enable gradient-based design optimization, the propagation of uncertainties, and the end-to-end training of scientific machine learning systems via backpropagation. 
For smooth flows, this can be carried out by either deriving a partial differential equation for the sensitivities (``differentiate then discretize'') or by automatic differentiation through a discrete numerical algorithm (``discretize then differentiate'').
The treatment of shocks is significantly more complicated.
The entropy condition defining vanishing viscosity is a differential \emph{in}equality and thus not amenable to the implicit function theorem. 
As such, formal differentiation through the equation does not yield the correct sensitivities. 
Likewise, differentiating through limiters results in spurious sensitivities \cite{lozano2019watch,bodony2022adjoint}, leading many practitioners to accept the systematic error incurred by ``freezing'' the limiters during differentiation \cite{nemec2006aerodynamic}.
PDE-based viscous regularizations avoid this problem \cite{bodony2022adjoint}, but frequently lack robustness at high Mach numbers \cite{glaubitz2019smooth}.
These challenges require delicate numerical and analytical treatment \cite{lozano2016note,murthy2026resolvent,ulbrich2002sensitivity,ulbrich2003adjoint,giles1997adjoint,giles1998properties}.

\subsection{Information geometric regularization (IGR)}
The limitations of existing viscous regularizations motivated the search for alternative regularization mechanisms. 
First attempts based on Leray regularization \cite{leray1934essai} failed to produce the correct shock speeds \cite{bhat2009regularization}. 
Subsequent attempts based on Hamiltonian mechanics \cite{guelmame2022hamiltonian} were insufficient to prevent singularity formation \cite{pu2018weakly,liu2019well,guelmame2022global}.
Information geometric regularization (IGR) overcomes these challenges, providing the first \emph{inviscid} regularization of compressible fluid mechanics \cite{cao2023information}.
It complements the physical pressure $P$ with an \emph{entropic pressure} $\Sigma$ defined in terms of an auxiliary elliptic equation, yielding 
\begin{equation}
\label{eqn:igr_euler}
\begin{cases}
     \partial_t \vct{\mu} + \divergence \left(\vct{\mu} \vct{u}^{\tp} + (P + \Sigma) \Id \right) &= \vct{f}_{\vct{\mu}}, \\
     \partial_t \rho + \divergence\left( \rho \vct{u}^{\tp} \right) &= f_\rho, \\
     \partial_t E + \divergence ((E + P + \Sigma)\vct{u}^{\tp}) &= \vct{u}^{\tp} \vct{f}_{\vct{\mu}} + f_E, \\
     \frac{\Sigma}{\rho} - \alpha \divergence \left(\frac{\nabla \Sigma}{\rho} \right) &= \alpha \left(\trace^2\left(\cnst{D}\vct{u}\right) + \trace\left(\left(\cnst{D}\vct{u}\right)^2\right)  \right).
\end{cases}
\end{equation}
Solutions of the IGR equations replace shocks with smooth profiles of width proportional to $\sqrt{\alpha}$, without damping oscillatory features \cite{cao2023information,barham2025hamiltonian}.
In the one-dimensional pressureless case, Cao and Sch\"afer \cite{cao2024information} show that IGR has global strong solutions that converge, in the limit of vanishing $\alpha$, to entropy solutions of the nominal system.
Barham et al.\ \cite{barham2026shock} prove the existence of traveling-wave solutions for a wide range of equations of state.
Eyob et al.\ \cite{eyob2026discontinuous} and Radhakrishnan et al.\ \cite{radhakrishnan2026shocks} introduce discontinuous Galerkin and finite-volume discretizations of IGR, with the latter enabling the first-ever numerical simulation exceeding a quadrillion degrees of freedom \cite{wilfong2025simulating}.
Recent works explore Hamiltonian \cite{barham2025hamiltonian} and thermodynamic \cite{taylor2026thermodynamically} extensions of IGR, as well as directional refinements \cite{xu2026compression}.

\subsection{Contribution in this work: (Adjoint) sensitivities of IGR}
Even in high-Mach-number flows, IGR maintains smooth solutions without damping out oscillatory features.
This makes it a promising approach to computing accurate sensitivities of flows with shocks. 
The present work realizes this potential through two main contributions. 
First, it derives the PDEs governing forward and adjoint sensitivities of IGR for periodic boundary conditions. 
Second, it demonstrates that these PDE-based sensitivities are consistent with those obtained by finite differences and automatic differentiation through a discontinuous Galerkin discretization.

\section{Forward Sensitivities}
\label{sec:forward_sensitivities}

\subsection{Overview}

We begin by deriving the so-called forward sensitivities of \cref{eqn:igr_euler}. 
Consider a single parameter $\theta$ on which the solution depends, for instance through its initial conditions, the forcing terms, or the equation of state. 
The forward sensitivities are derivatives of the solutions $\vct{\mu}$, $\rho$, and $E$ with respect to $\theta$. 
For any variable, such as $\vct{\mu}$ or $P$, we denote as $\sens{\vct{\mu}}$ or $\sens{P}$ the corresponding sensitivity variable $\partial_\theta \vct{\mu}$ or $\partial_\theta P$.

\subsection{Derivation}
We derive the forward sensitivity equations by differentiating \cref{eqn:igr_euler} with respect to $\theta$.
The derivative of the momentum balance equation yields
\begin{equation}
     \partial_t \sens{\vct{\mu}} + \divergence \left(\sens{\vct{\mu}} \vct{u}^{\tp} + \vct{\mu} \sens{\vct{u}}^{\tp} + (\sens{P} + \sens{\Sigma}) \Id \right) = \sens{\vct{f}}_{\vct{\mu}}.
\end{equation}
Likewise, the mass balance equation yields
\begin{equation}
     \partial_t \sens{\rho} + \divergence\left( \sens{\rho} \vct{u}^{\tp} + \rho \sens{\vct{u}}^{\tp} \right) = \sens{f}_\rho,
     \quad \text{which simplifies to} \quad 
     \partial_t \sens{\rho} + \divergence\left( \sens{\vct{\mu}}^{\tp} \right) = \sens{f}_\rho.
\end{equation}
Finally, differentiating the energy equation gives 
\begin{equation}
     \partial_t \sens{E} + \divergence ((\sens{E} + \sens{P} + \sens{\Sigma})\vct{u}^{\tp} + (E + P + \Sigma) \sens{\vct{u}}^{\tp}) = \sens{\vct{u}}^{\tp} \vct{f}_{\vct{\mu}} + \vct{u}^{\tp} \sens{\vct{f}}_{\vct{\mu}} + \sens{f}_E.
\end{equation}
The above yields an evolution PDE for the sensitivities $\sens{\vct{\mu}}$, $\sens{\rho}$, and $\sens{E}$. 
But besides the forward solution variables, the sensitivity equations also depend on the sensitivities $\sens{\vct{u}}$, $\sens{P}$, and $\sens{\Sigma}$ of the velocity, pressure, and entropic pressure.
We can compute these sensitivities from the primal solution and the sensitivity variables $\sens{\vct{\mu}}$, $\sens{\rho}$, and $\sens{E}$ as follows.
For the velocity sensitivity, we use the quotient rule to compute 
\begin{equation}
\sens{\vct{u}} = \partial_\theta \left( \frac{\vct{\mu}}{\rho} \right) = \frac{\sens{\vct{\mu}}}{\rho} - \frac{\vct{\mu} \sens{\rho}}{\rho^2}.
\end{equation}
For the pressure sensitivity $\sens{P}$, we apply the chain rule to the equation of state $P = P(\rho, e)$, where $e = E/\rho - \vct{u}^{\tp}\vct{u}/2$ is the specific internal energy.
Often, $P$ depends only on the density and internal energy (and possibly directly on $\theta$), yielding
\begin{equation}
\sens{P} = \frac{\partial P}{\partial \rho} \sens{\rho} + \frac{\partial P}{\partial e} \sens{e} + \frac{\partial P}{\partial \theta},
\quad \text{where} \quad
\sens{e} = \partial_\theta \left( \frac{E}{\rho} - \frac{\vct{u}^{\tp} \vct{u}}{2} \right) = \frac{\sens{E}}{\rho} - \frac{\sens{\rho} E}{\rho^2} - \vct{u}^{\tp} \sens{\vct{u}}.
\end{equation}
In polytropic equations of state $P = (\gamma - 1) \rho e$, the pressure sensitivity simplifies to
\begin{equation}\label{eqn:polytropic_pressure_sensitivity}
\sens{P} = (\gamma - 1) \left( \sens{\rho} e + \rho \sens{e} \right) + \sens{\gamma} \rho e.
\end{equation}
Finally, the entropic pressure sensitivity $\sens{\Sigma}$ solves
\begin{equation}
     \frac{\sens{\Sigma}}{\rho} - \alpha \divergence \left(\frac{\nabla \sens{\Sigma}}{\rho}\right) 
     = 2\alpha \left(\trace\left( \cnst{D}\vct{u}\right) \trace\left( \cnst{D}\sens{\vct{u}}\right) + \trace\left(\cnst{D}\vct{u} \cnst{D}\sens{\vct{u}}\right)\right) + \left(\frac{\sens{\rho}\Sigma}{\rho^2} - \alpha \divergence\left(\frac{\sens{\rho}\nabla \Sigma}{\rho^2}\right) \right).
\end{equation}
This is an elliptic PDE for $\sens{\Sigma}$ with the same structure as the primal $\Sigma$ equation, but with a different right-hand side.
Thus, it can reuse the elliptic solver used for the original IGR calculation.
In order to simplify the computation of the right-hand side (which currently requires second derivatives), we can apply the product rule to obtain
\begin{equation}
     \divergence\left(\frac{\sens{\rho}}{\rho} \frac{\nabla \Sigma}{\rho}\right) =  
     \nabla\left(\frac{\sens{\rho}}{\rho} \right) \cdot \frac{\nabla \Sigma}{\rho} + \frac{\sens{\rho}}{\rho} \divergence\left(\frac{\nabla \Sigma}{\rho}\right).
\end{equation}
By plugging in the definition of $\Sigma$, this yields
\begin{align}
     \frac{\sens{\Sigma}}{\rho} - \alpha \divergence \left(\frac{\nabla \sens{\Sigma}}{\rho}\right)
     & = \alpha \Bigg(\trace\left( \cnst{D}\vct{u}\right) \trace\left(2 \cnst{D}\sens{\vct{u}} 
     +  \frac{\sens{\rho}}{\rho} \cnst{D} \vct{u}\right) \\ 
     & \qquad + \trace\left(\cnst{D}\vct{u} \left(2 \cnst{D}\sens{\vct{u}} + \frac{\sens{\rho}}{\rho}\cnst{D}\vct{u}\right)\right) - \nabla\left(\frac{\sens{\rho}}{\rho}\right) \cdot \frac{\nabla \Sigma}{\rho} \Bigg).
\end{align}
The combined forward sensitivity PDE reads
{\small
\begin{equation}
\label{eqn:forward_sensitivity_system}
\begin{cases}
     \partial_t \vct{\mu} + \divergence \left(\vct{\mu} \vct{u}^{\tp} + (P + \Sigma) \Id \right) &= \vct{f}_{\vct{\mu}}, \\
     \partial_t \rho + \divergence\left( \rho \vct{u}^{\tp} \right) &= f_\rho, \\
     \partial_t E + \divergence ((E + P + \Sigma)\vct{u}^{\tp}) &= \vct{u}^{\tp} \vct{f}_{\vct{\mu}} + f_E, \\
     \partial_t \sens{\vct{\mu}} + \divergence \left(\sens{\vct{\mu}} \vct{u}^{\tp} + \vct{\mu} \sens{\vct{u}}^{\tp} + (\sens{P} + \sens{\Sigma}) \Id \right) &= \sens{\vct{f}}_{\vct{\mu}},\\
     \partial_t \sens{\rho} + \divergence\left( \sens{\rho} \vct{u}^{\tp} + \rho \sens{\vct{u}}^{\tp}\right) &= \sens{f}_\rho,\\
     \partial_t\sens{E} + \divergence ((E + P + \Sigma)\sens{\vct{u}}^{\tp} + (\sens{E} + \sens{P} + \sens{\Sigma})\vct{u}^{\tp}) &= \sens{\vct{u}}^{\tp} \vct{f}_{\vct{\mu}} + \vct{u}^{\tp} \sens{\vct{f}}_{\vct{\mu}} + \sens{f}_E,\\
     \frac{\Sigma}{\rho} - \alpha \divergence \left(\frac{\nabla \Sigma}{\rho} \right) &= \alpha \left(\trace^2\left(\cnst{D}\vct{u}\right) + \trace\left(\left(\cnst{D}\vct{u}\right)^2\right)\right),\\
     \frac{\sens{\Sigma}}{\rho} - \alpha \divergence \left(\frac{\nabla \sens{\Sigma}}{\rho}\right) 
     &= \alpha \Bigg(\trace\left( \cnst{D}\vct{u}\right) \trace\left(2 \cnst{D}\sens{\vct{u}} + \frac{\sens{\rho}}{\rho}\cnst{D} \vct{u}\right) \\  &+\! \trace\left(\cnst{D}\vct{u} \left(2 \cnst{D}\sens{\vct{u}} \!+\! \frac{\sens{\rho}}{\rho} \cnst{D}\vct{u}\right)\right) \!-\! \nabla\left(\frac{\sens{\rho}}{\rho}\right) \!\cdot\! \frac{\nabla \Sigma}{\rho} \Bigg).
\end{cases}
\end{equation}
}

\section{Adjoint Sensitivities}
\label{sec:adjoint_sensitivities}

\subsection{Overview}
The forward sensitivities derived above provide the rate of change of solutions under a single parameter perturbation. 
But many applications require sensitivities of a scalar (or low-dimensional) functional of the solution with respect to (possibly infinitely) many parameters.
Adjoint methods avoid the cost of computing sensitivities for each parameter by deriving a partial differential equation for the gradient of an output functional. 
We consider linear functionals of the form
\begin{equation}
     J\left(
          \begin{pmatrix}
          \vct{\mu} \\ 
          \rho \\
          E 
          \end{pmatrix}
     \right)
          = \int \limits_0^T \int \limits_\Omega 
          \begin{pmatrix}
               \vct{g}_{\vct{\mu}}\left(\vct{x}, t\right) \\ 
               g_\rho\left(\vct{x}, t\right) \\
               g_E\left(\vct{x}, t\right)
          \end{pmatrix}
          \cdot 
          \begin{pmatrix}
               \vct{\mu}\left(\vct{x}, t\right) \\
               \rho\left(\vct{x}, t\right) \\
               E\left(\vct{x}, t\right)
          \end{pmatrix}
          \D \vct{x} \D t
          + \int \limits_\Omega 
          \begin{pmatrix}
          \vct{h}_{\vct{\mu}}\left(\vct{x}\right) \\ 
          h_\rho\left(\vct{x}\right) \\
          h_E\left(\vct{x}\right)
          \end{pmatrix}
          \cdot 
          \begin{pmatrix}
          \vct{\mu}\left(\vct{x}, T\right) \\
          \rho\left(\vct{x}, T\right) \\
          E\left(\vct{x}, T\right)
          \end{pmatrix}
          \D \vct{x}
\end{equation}
that combine a running loss (integrated over time) with a terminal loss (evaluated at the final time $T$).
Nonlinear functionals can be treated by linearization.
The goal of the adjoint state method is to derive a PDE for the adjoint variables $\adj{\vct{\mu}}$, $\adj{\rho}$, and $\adj{E}$ such that the derivative of $J$ with respect to any perturbation of $\theta$ can be computed as a simple inner product of the adjoint variables with the perturbation.
More precisely, we consider the IGR equations with parametric right-hand sides, of the form
\begin{equation}
\label{eqn:igr_euler_parametric}
\begin{cases}
     \partial_t \vct{\mu} + \divergence \left(\vct{\mu} \vct{u}^{\tp} + (P + \Sigma) \Id \right) &= \vct{f}_{\vct{\mu}} + \vct{\theta}_{\vct{\mu}}, \\
     \partial_t \rho + \divergence\left( \rho \vct{u}^{\tp} \right) &= f_\rho + \theta_\rho, \\
     \partial_t E + \divergence ((E + P + \Sigma)\vct{u}^{\tp}) &= \vct{u}^{\tp} \vct{f}_{\vct{\mu}} + f_E + \theta_E, \\
     \frac{\Sigma}{\rho} - \alpha \divergence \left(\frac{\nabla \Sigma}{\rho} \right) &= \alpha \left(\trace^2\left(\cnst{D}\vct{u}\right) + \trace\left(\left(\cnst{D}\vct{u}\right)^2\right)  \right) 
\end{cases}
\end{equation}
on a periodic domain $\Omega$ and initial conditions of the form 
\begin{equation}
     \vct{\mu}(\vct{x}, 0) = \vct{\mu}_0(\vct{x}) + \vct{\theta}_{\vct{\mu}, 0}(\vct{x}), \quad
     \rho(\vct{x}, 0) = \rho_0(\vct{x}) + \theta_{\rho, 0}(\vct{x}), \quad
     E(\vct{x}, 0) = E_0(\vct{x}) + \theta_{E, 0}(\vct{x}).
\end{equation}
We then seek a PDE that governs the evolution of the adjoint variables $\adj{\vct{\mu}}$, $\adj{\rho}$, and $\adj{E}$, allowing us to compute the functional derivative of $J$ with respect to any variation as
\begin{equation}
     \frac{\D J}{\D \theta} = \int \limits_0^T \int \limits_\Omega 
     \begin{pmatrix}
          \adj{\vct{\mu}}\left(\vct{x}, t\right) \\ 
          \adj{\rho}\left(\vct{x}, t\right) \\
          \adj{E}\left(\vct{x}, t\right)
     \end{pmatrix}
     \cdot 
     \begin{pmatrix}
          \vct{\theta}_{\vct{\mu}}\left(\vct{x}, t\right) \\
          \theta_\rho\left(\vct{x}, t\right) \\
          \theta_E\left(\vct{x}, t\right)
     \end{pmatrix}
     \D \vct{x} \D t
     + 
     \int \limits_\Omega
     \begin{pmatrix}
          \adj{\vct{\mu}}\left(\vct{x}, 0\right) \\ 
          \adj{\rho}\left(\vct{x}, 0\right) \\
          \adj{E}\left(\vct{x}, 0\right)
     \end{pmatrix}
     \cdot 
     \begin{pmatrix}
          \vct{\theta}_{\vct{\mu}, 0}\left(\vct{x}\right) \\
          \theta_{\rho, 0}\left(\vct{x}\right) \\
          \theta_{E, 0}\left(\vct{x}\right)
     \end{pmatrix}
     \D \vct{x}.
\end{equation}
This allows us to compute the derivative of $J$ with respect to any parameter $\theta$ by integrating the adjoint variables against the parameter perturbation, without needing to compute the forward sensitivities for each parameter separately.

\subsection{Adjoint derivation in the abstract, for local fluxes}
We first carry out the adjoint derivation abstractly, for a hyperbolic conservation law of the form
\begin{equation}
     \partial_t \vct{q} + \divergence \left(\vct{F}\left(\vct{q}\right) \right) = \vct{f} + \vct{\theta}, \quad \vct{q}(\vct{x}, 0) = \vct{q}_0(\vct{x}) + \vct{\theta}_0(\vct{x})
\end{equation}
with $k$ components and in $d$ spatial dimensions, and an output functional of the form
\begin{equation}
     J\left(\vct{q}\right) = \int \limits_0^T \int \limits_\Omega \vct{g}\left(\vct{x}, t\right) \cdot \vct{q}\left(\vct{x}, t\right) \D \vct{x} \D t + \int \limits_\Omega \vct{h}\left(\vct{x}\right) \cdot \vct{q}\left(\vct{x}, T\right) \D \vct{x}.
\end{equation}
Here, we take $\divergence$ to be the row-wise divergence operator and $\vct{F}: \R^k \to \R^{k \times d}$ to be the flux function.
To simplify the calculations, we assume that the initial conditions are not parametrized, i.e., $\vct{\theta}_0 = 0$, and that the output functional does not depend on the terminal state, i.e., $\vct{h} = 0$.
The adjoint state variable $\adj{\vct{q}}$ then satisfies
\begin{equation}
     \frac{\D J}{\D \theta}\left(\delta \theta\right) 
     = \int \limits_0^T \int \limits_\Omega \adj{\vct{q}}\left(\vct{x}, t\right) \cdot \delta \theta\left(\vct{x}, t\right) \D \vct{x} \D t
\end{equation}
At the same time, the variation of $J$ with respect to $\theta$ can be computed by first computing the forward sensitivity $\sens{\vct{q}}$ of the solution $\vct{q}$ with respect to $\theta$, and then integrating against the output functional gradient $\vct{g}$, yielding
\begin{equation}
     \frac{\D J}{\D \theta}\left(\delta \theta\right) = \int \limits_0^T \int \limits_\Omega \vct{g}\left(\vct{x}, t\right) \cdot \sens{\vct{q}}\left(\vct{x}, t\right) \D \vct{x} \D t.
\end{equation}
The forward sensitivity, in turn, satisfies the PDE 
\begin{equation}
     \partial_t \sens{\vct{q}} + \divergence \left( \cnst{D}_{\vct{q}}\vct{F}\left(\vct{q}\right) \sens{\vct{q}} \right) = \delta \theta, \quad \sens{\vct{q}}(\vct{x}, 0) = 0.
\end{equation}
We can plug the above in for the $\delta \theta$ in the expression for $\D J/\D \theta$ to obtain
\begin{equation}
     \int \limits_0^T \int \limits_\Omega \vct{g} \cdot \sens{\vct{q}} \D \vct{x} \D t   
     = \int \limits_0^T \int \limits_\Omega \adj{\vct{q}} \cdot \left(\partial_t \sens{\vct{q}} + \divergence \left( \cnst{D}_{\vct{q}}\vct{F}\left(\vct{q}\right) \sens{\vct{q}} \right)\right) \D \vct{x} \D t
\end{equation}
We now integrate by parts in time and space to move the derivatives from $\sens{\vct{q}}$ onto $\adj{\vct{q}}$,
\begin{align}
     &\int \limits_0^T \int \limits_\Omega \vct{g} \cdot \sens{\vct{q}} \D \vct{x} \D t 
     = \int \limits_\Omega \int \limits_0^T  - \partial_t \adj{\vct{q}} \cdot \sens{\vct{q}} - \cnst{D} \adj{\vct{q}} \mathrel{\colon} \left(\cnst{D}_{\vct{q}}\vct{F}\left(\vct{q}\right) \cdot \sens{\vct{q}} \right)\D t 
          + \adj{\vct{q}}\left(\vct{x}, T\right) \cdot \sens{\vct{q}}\left(\vct{x}, T\right) \D \vct{x} \\
     &\quad = \int \limits_0^T \int \limits_\Omega \left(- \partial_t \adj{\vct{q}}\left(\vct{x}, t\right) - \nabla_{\vct{q}}\left(  \vct{F}\left(\vct{q}\right) \mathrel{\colon} \cnst{D} \adj{\vct{q}} \right)\right) \cdot \sens{\vct{q}} \D \vct{x} \D t
     + \int \limits_\Omega \adj{\vct{q}}\left(\vct{x}, T\right) \cdot \sens{\vct{q}}\left(\vct{x}, T\right) \D \vct{x},
\end{align}
using periodicity to ignore boundary terms.
The main challenge in deriving the above expression is to keep track of the contraction involving the third-order tensor $\cnst{D}_{\vct{q}}\vct{F}(\vct{q})$.
The above expression can be deduced by checking dimensional consistency and remembering that $\sens{\vct{q}}$ is contracted over the derivatives with respect to $\vct{q}$.
Since $\vct{\theta}$ enters as a source term in the evolution equation, any possible variation can at least formally be realized by some $\delta \vct{\theta}$.
Thus, requiring the identity to hold for all $\sens{\vct{q}}$ yields
\begin{equation}
     \label{eqn:adjoint_abstract_local_nonconservative}
     - \partial_t \adj{\vct{q}} - \nabla_{\vct{q}}\left(  \vct{F}\left(\vct{q}\right) \mathrel{\colon} \cnst{D} \adj{\vct{q}} \right) = \vct{g}, \quad \adj{\vct{q}}\left(\vct{x}, T\right) \equiv 0.
\end{equation}
As is typical for adjoint state equations of nonlinear conservation laws, the adjoint PDE is not in conservation form. 
We can instead formulate it as a conservation law with a source term, by applying the product rule to the flux term, yielding
\begin{equation}
     \label{eqn:adjoint_abstract_local_conservative}
     - \partial_t \adj{\vct{q}} - \divergence\left( \nabla_{\vct{q}} \left(\adj{\vct{q}}^\tp \vct{F}\left(\vct{q}\right) \right) \right) + \left(\adj{\vct{q}}^{\tp} \divergence\left(\cnst{D}_{\vct{q}}\vct{F}\left(\vct{q}\right)\right)\right) = \vct{g}, \quad \adj{\vct{q}}\left(\vct{x}, T\right) \equiv 0.
\end{equation}
The limit $\vct{g}(\vct{x}, t) \rightarrow \delta_{T}(t) \vct{h}(\vct{x})$ implies that in general, $\adj{\vct{q}}\left(\vct{x}, T\right) \equiv \vct{h}(\vct{x})$.

\subsection{Extension to nonlocal fluxes}
The above derivation assumed that the flux $\vct{F}$ is a local function of the state $\vct{q}$, i.e., that $\vct{F}(\vct{q})(\vct{x})$ depends on $\vct{q}$ only through its value $\vct{q}(\vct{x})$ at $\vct{x}$.
This holds for most commonly considered hyperbolic conservation laws.
But as discussed in \cite[Section 7.1]{cao2023information}, the nonlocal dependence of the flux is a crucial feature of IGR.
Since the nonlocality of the IGR flux is mediated by only $\Sigma$, we can separate the treatment of the nonlocal and local components by considering an $\vct{F}$ that depends locally on both $\vct{q}$ and $\Sigma$, with a nonlocal dependence of $\Sigma$ on $\vct{q}$.
For a fixed field $\vct{x} \mapsto \underline{\Sigma}(\vct{x})$, the map $\vct{q} \mapsto \vct{F}(\vct{q}, \underline{\Sigma})$ is an ordinary local flux, while for a fixed $\vct{x} \mapsto \underline{\vct{q}}(\vct{x})$, the map $\vct{q} \mapsto \vct{F}(\underline{\vct{q}}, \Sigma(\vct{q}))$ carries the entire nonlocal dependence.
We define gradients of these restricted maps, evaluated at the actual solution, as
\begin{align}
     \nabla^{\Sigma}_{\vct{q}} \left(\vct{F} \mathrel{\colon} \cnst{D}\adj{\vct{q}}\right)
     &\defeq \left.\nabla_{\vct{q}} \left(\vct{F}\left(\vct{q}, \underline{\Sigma}\right) \mathrel{\colon} \cnst{D}\adj{\vct{q}}\right)\right|_{\underline{\Sigma} = \Sigma\left(\vct{q}\right)},\\
     \nabla^{\vct{q}}_{\vct{q}} \left(\vct{F} \mathrel{\colon} \cnst{D}\adj{\vct{q}}\right)
     &\defeq \left.\nabla_{\vct{q}} \left(\vct{F}\left(\underline{\vct{q}}, \Sigma\left(\vct{q}\right)\right) \mathrel{\colon} \cnst{D}\adj{\vct{q}}\right)\right|_{\underline{\vct{q}} = \vct{q}},
\end{align}
where the superscript indicates the quantity held fixed under differentiation.
The first gradient is a pointwise gradient, while the second is the $L^2$ gradient of $\vct{q} \mapsto \int_{\Omega} \vct{F}\left(\underline{\vct{q}}, \Sigma(\vct{q})\right) \mathrel{\colon} \cnst{D} \adj{\vct{q}} \D \vct{x}$.
By the chain rule, $\nabla_{\vct{q}} = \nabla^{\Sigma}_{\vct{q}} + \nabla^{\vct{q}}_{\vct{q}}$, so that \cref{eqn:adjoint_abstract_local_nonconservative} becomes
\begin{equation}
     \label{eqn:adjoint_abstract_nonlocal_nonconservative}
     - \partial_t \adj{\vct{q}} - \nabla^{\Sigma}_{\vct{q}} \left(\vct{F} \mathrel{\colon} \cnst{D} \adj{\vct{q}} \right)  - \nabla^{\vct{q}}_{\vct{q}} \left(\vct{F} \mathrel{\colon} \cnst{D} \adj{\vct{q}} \right)= \vct{g},
\end{equation}
with the same terminal condition as before.
Since $\Sigma$ is held fixed in the first term, the passage from \cref{eqn:adjoint_abstract_local_nonconservative} to \cref{eqn:adjoint_abstract_local_conservative} applies to it verbatim, yielding the conservative form
\begin{equation}
     \label{eqn:adjoint_abstract_nonlocal_conservative}
     - \partial_t \adj{\vct{q}} - \divergence\left( \nabla^{\Sigma}_{\vct{q}} \left(\adj{\vct{q}}^\tp \vct{F} \right) \right) + \adj{\vct{q}}^{\tp} \divergence\left(\cnst{D}^{\Sigma}_{\vct{q}}\vct{F}\right) - \nabla^{\vct{q}}_{\vct{q}} \left(\vct{F} \mathrel{\colon} \cnst{D} \adj{\vct{q}} \right)= \vct{g},
\end{equation}
with the gradient $\nabla^{\Sigma}_{\vct{q}}\left(\adj{\vct{q}}^\tp \vct{F}\right)$ and the Jacobian $\cnst{D}^{\Sigma}_{\vct{q}}\vct{F}$ at frozen $\Sigma$ defined analogously.
\subsection{Derivation of the IGR adjoint equations in advective form}
We now apply the above abstract derivation to the information-geometrically regularized Euler equations.
We avoid introducing third-order tensors by always first contracting the flux with the adjoint variable $\adj{\vct{q}}$ and then computing the gradient with respect to the state variable $\vct{q}$.
We write $\bar{P} = P + \Sigma$ and observe that the state variable is 
\begin{equation}
     \vct{q} = 
     \begin{pmatrix}
     \vct{\mu} \\ 
     \rho\\
     E
     \end{pmatrix}
     \quad \text{and the flux function} \quad
     \vct{F}(\vct{q}) = 
     \begin{pmatrix}
          \vct{\mu} \vct{u}^{\tp} + \bar{P} \Id \\
          \rho \vct{u}^{\tp} \\
          (E + \bar{P})\vct{u}^{\tp}
     \end{pmatrix}.
\end{equation}
We first compute the advective form of the adjoint PDE.
We use subscripts to denote gradients and derivatives with respect to state variables, as in $P_{\vct{\mu}} = \nabla_{\vct{\mu}} P$, $P_\rho = \partial_\rho P$, $P_E = \partial_E P$, and define the generalized enthalpy as $\bar{H} = E + P + \Sigma$.
Since the elliptic equation in \cref{eqn:igr_euler} does not involve $E$, we have $\Sigma_E = 0$.
We begin with the first term of \cref{eqn:adjoint_abstract_nonlocal_nonconservative}, the gradient at frozen $\Sigma$.
Since $\Sigma$ is frozen, the derivatives of $\bar{P}$ and $\bar{H}$ reduce to those of $P$, i.e., $\bar{P}_{\vct{\mu}} = P_{\vct{\mu}}$, $\bar{P}_\rho = \bar{H}_\rho = P_\rho$, and $\bar{P}_E = P_E$, yielding
\begin{align}
     &\nabla_{\vct{q}}^{\Sigma}\left(
          \int \limits_\Omega
          \begin{pmatrix}
               \vct{\mu} \vct{u}^{\tp} + \bar{P} \Id \\
               \rho \vct{u}^{\tp} \\
               (E + \bar{P})\vct{u}^{\tp}
          \end{pmatrix}
          \mathrel{\colon} 
          \begin{pmatrix}
               \cnst{D} \adj{\vct{\mu}} \\
               \cnst{D} \adj{\rho} \\
               \cnst{D} \adj{E}
          \end{pmatrix}
          \D \vct{x}
     \right)\\
     =&
     \nabla_{\vct{q}}^{\Sigma}\left(
         \int \limits_\Omega
          \frac{\vct{\mu}^{\tp} \left[\cnst{D} \adj{\vct{\mu}}\right] \vct{\mu}}{\rho} + \bar{P} \divergence\left(\adj{\vct{\mu}}\right) + \left[\cnst{D}\adj{\rho}\right] \vct{\mu} + \left[\cnst{D} \adj{E}\right] \vct{\mu} \frac{\bar{H}}{\rho}
         \D \vct{x} 
     \right)\\
     =& 
     \begin{pmatrix}
          \left(\left[\cnst{D} \adj{\vct{\mu}}\right] + \left[\cnst{D} \adj{\vct{\mu}}\right]^{\tp}\right) \vct{u}  + \divergence\left(\adj{\vct{\mu}}\right) P_{\vct{\mu}}   + \nabla \adj{\rho} + \frac{\bar{H}}{\rho} \nabla \adj{E} + \frac{\left[\cnst{D}\adj{E}\right] \vct{\mu}}{\rho} P_{\vct{\mu}} \\
          - \frac{\vct{\mu}^{\tp} \left[\cnst{D} \adj{\vct{\mu}}\right]\vct{\mu}}{\rho^2}  + \left(P_\rho\right) \divergence\left(\adj{\vct{\mu}}\right)
          + \left[\cnst{D}\adj{E}\right] \vct{\mu} \left(\frac{P_\rho}{\rho}
          - \frac{\bar{H}}{\rho^2}\right)  \\
          P_E \left( \divergence\left(\adj{\vct{\mu}}\right) + \frac{\left[\cnst{D}\adj{E}\right] \vct{\mu}}{\rho}\right) + \frac{\left[\cnst{D}\adj{E}\right] \vct{\mu}}{\rho}
     \end{pmatrix}.
\end{align}
We now compute the second term by freezing all $\vct{q}$ except those entering as the argument of $\Sigma$, denoted by $\nabla_{\vct{q}}^{\vct{q}}$, obtaining 
\begin{equation}
     \nabla_{\vct{q}}^{\vct{q}}\!\left(
         \int \limits_\Omega
          \frac{\vct{\mu}^{\tp} \left[\cnst{D} \adj{\vct{\mu}}\right] \vct{\mu}}{\rho} + \bar{P} \divergence\left(\adj{\vct{\mu}}\right) + \left[\cnst{D}\adj{\rho}\right] \vct{\mu} + \left[\cnst{D} \adj{E}\right] \vct{\mu} \frac{\bar{H}}{\rho}
         \D \vct{x} 
     \!\right)
     \!=\!
     \left(\!\divergence\left(\adj{\vct{\mu}}\right) + \frac{\left[\cnst{D} \adj{E}\right]\vct{\mu}}{\rho}\!\right)
     \!\Sigma_{\vct{q}},
\end{equation}
where the right hand side is interpreted as the $L^2$ representative of 
\begin{equation}
     \sens{\vct{q}} \mapsto \left. \frac{\D}{\D \delta} \int \limits_\Omega \left(\!\divergence\left(\adj{\vct{\mu}}\right) + \frac{\left[\cnst{D} \adj{E}\right]\vct{\mu}}{\rho}\!\right)
     \!\Sigma(\vct{q} + \delta \sens{\vct{q}}) \D \vct{x} \right|_{\delta = 0}.
\end{equation}
Together, we obtain \cref{eqn:adjoint_abstract_nonlocal_nonconservative} in the form
\begin{equation}
     \begin{split}
          - 
          \partial_t 
          \begin{pmatrix}
               \adj{\vct{\mu}}\\
               \adj{\rho}\\
               \adj{E}
          \end{pmatrix}
          - 
          \begin{pmatrix}
               \left(\left[\cnst{D} \adj{\vct{\mu}}\right] + \left[\cnst{D} \adj{\vct{\mu}}\right]^{\tp}\right) \vct{u}  
                    + \divergence\left(\adj{\vct{\mu}}\right) P_{\vct{\mu}}
                    + \nabla \adj{\rho} + \frac{\bar{H}}{\rho} \nabla \adj{E}
                    + \frac{\left[\cnst{D}\adj{E}\right] \vct{\mu}}{\rho}P_{\vct{\mu}}\\
               - \frac{\vct{\mu}^{\tp} \left[\cnst{D} \adj{\vct{\mu}}\right]\vct{\mu}}{\rho^2}  
                    + \left(P_\rho\right) \divergence\left(\adj{\vct{\mu}}\right) 
                    + \left[\cnst{D}\adj{E}\right] \vct{\mu} \left(\frac{P_\rho}{\rho} - \frac{\bar{H}}{\rho^2}\right)  \\
               P_E \left( \divergence\left(\adj{\vct{\mu}}\right) + \frac{\left[\cnst{D}\adj{E}\right] \vct{\mu}}{\rho}\right) + \frac{\left[\cnst{D}\adj{E}\right] \vct{\mu}}{\rho}
          \end{pmatrix} \qquad \qquad\\
          -
          \left(\divergence\left(\adj{\vct{\mu}}\right) + \frac{\left[\cnst{D} \adj{E}\right]\vct{\mu}}{\rho}\right)
          \Sigma_{\vct{q}} = 
          \begin{pmatrix}
               \vct{g}_{\vct{\mu}}\left(\vct{x}, t\right) \\ 
               g_\rho\left(\vct{x}, t\right) \\
               g_E\left(\vct{x}, t\right)
          \end{pmatrix}
          .
     \end{split}
\end{equation}
\subsection{Derivation of the IGR adjoint equations in conservative form}
Now we move on to the conservative form of the adjoint PDE, which is more suitable for integration into finite-volume solvers.
We begin by computing the flux term  
\begin{align}
&\nabla_{\vct{q}}^{\Sigma} \left( 
\begin{pmatrix}
     \adj{\vct{\mu}}^\tp &
     \adj{\rho} &
     \adj{E}
\end{pmatrix}
\begin{pmatrix}
     \vct{\mu} \vct{u}^{\tp} + \bar{P} \Id \\
     \rho \vct{u}^{\tp} \\
     (E + \bar{P})\vct{u}^{\tp}
\end{pmatrix}
\right) =
\nabla_{\vct{q}}^{\Sigma} 
\left(
\adj{\vct{\mu}}^\tp \vct{\mu} \frac{\vct{\mu}^{\tp}}{\rho} + \bar{P} \adj{\vct{\mu}}^{\tp} + \adj{\rho} \vct{\mu}^{\tp} + \adj{E} \frac{\bar{H}\vct{\mu}^{\tp}}{\rho} 
\right)\\
=& 
\begin{pmatrix}
     \frac{\adj{\vct{\mu}} \vct{\mu}^\tp + \adj{\vct{\mu}}^{\tp} \vct{\mu} \Id}{\rho} + \bar{P}_{\vct{\mu}} \adj{\vct{\mu}}^{\tp} + \adj{\rho} \Id + \adj{E} \frac{\bar{H}}{\rho} \Id + \adj{E} \frac{\bar{P}_{\vct{\mu}} \vct{\mu}^{\tp}}{\rho}\\
     - \frac{\adj{\vct{\mu}}^\tp \vct{\mu} \vct{\mu}^\tp}{\rho^2} + \bar{P}_{\rho} \adj{\vct{\mu}}^{\tp} + \adj{E} \left(\partial_\rho\left(\frac{\bar{H}}{\rho}\right)\right) \vct{\mu}^{\tp}\\
     \bar{P}_E \adj{\vct{\mu}}^{\tp} + \adj{E} \frac{\left(1 + \bar{P}_E\right) \vct{\mu}^{\tp}}{\rho}
\end{pmatrix}\\
=& 
\begin{pmatrix}
     \frac{\adj{\vct{\mu}} \vct{\mu}^\tp + \adj{\vct{\mu}}^{\tp} \vct{\mu} \Id}{\rho} 
          + P_{\vct{\mu}} \adj{\vct{\mu}}^{\tp}
          + \adj{\rho} \Id 
          + \adj{E} \frac{\bar{H}}{\rho} \Id 
          + \adj{E} \frac{P_{\vct{\mu}}\vct{\mu}^{\tp}}{\rho}\\
     - \frac{\adj{\vct{\mu}}^\tp \vct{\mu} \vct{\mu}^\tp}{\rho^2} 
          + P_{\rho} \adj{\vct{\mu}}^{\tp} 
          + \adj{E} \left(\left(\frac{P_{\rho}}{\rho}\right) - \left(\frac{\bar{H}}{\rho^2}\right)\right) \vct{\mu}^{\tp}\\
     P_E \adj{\vct{\mu}}^{\tp} + \adj{E} \frac{\left(1 + P_E \right) \vct{\mu}^{\tp}}{\rho}
\end{pmatrix}
\end{align}
We now compute the source term $\adj{\vct{q}}^{\tp}\divergence(\cnst{D}^{\Sigma}_{\vct{q}}\vct{F})$ of \cref{eqn:adjoint_abstract_nonlocal_conservative}. 
We avoid third-order tensors by first rewriting $\adj{\vct{q}}^{\tp}\divergence(\vct{F})$, denoting matrices as $(\cdots)$ and line breaks as $[\cdots]$,
\begin{align}
     & \begin{pmatrix}
               \adj{\vct{\mu}}^\tp &
               \adj{\rho} &
               \adj{E}
          \end{pmatrix}
          \!
          \begin{pmatrix}
               \divergence\left(\vct{\mu} \vct{u}^{\tp} + \bar{P} \Id \right) \\
               \divergence\left(\rho \vct{u}^{\tp} \right) \\
               \divergence\left((E + \bar{P})\vct{u}^{\tp} \right)
          \end{pmatrix}
     =
     \begin{pmatrix}
               \adj{\vct{\mu}}^\tp &
               \adj{\rho} &
               \adj{E}
          \end{pmatrix}
          \!
          \begin{pmatrix}
               \divergence\left(\frac{\vct{\mu} \vct{\mu}^{\tp}}{\rho} + \bar{P} \Id \right) \\
               \divergence\left(\vct{\mu} \right) \\
               \divergence\left(\frac{(E + \bar{P})\vct{\mu}^{\tp}}{\rho} \right)
          \end{pmatrix}\\
     &\qquad =
     \begin{pmatrix}
          \adj{\vct{\mu}}^\tp &
          \adj{\rho} &
          \adj{E}
     \end{pmatrix}
     \begin{pmatrix}
          \frac{\divergence(\vct{\mu})\vct{\mu} + (\vct{\mu}^{\tp} \nabla)\vct{\mu}}{\rho} - \frac{\vct{\mu} \vct{\mu}^{\tp}\nabla \rho}{\rho^2} + \nabla \bar{P}\\
          \divergence\left(\vct{\mu} \right) \\
          \frac{\vct{\mu}^{\tp} \nabla (E + \bar{P}) + (E + \bar{P}) \divergence(\vct{\mu})}{\rho} - \frac{(E + \bar{P})\vct{\mu}^{\tp} \nabla \rho}{\rho^2}
     \end{pmatrix}\\
     &\qquad=
     \begin{bmatrix}
           \frac{\adj{\vct{\mu}}^{\tp}\divergence(\vct{\mu})\vct{\mu} + \adj{\vct{\mu}}^{\tp}(\vct{\mu}^{\tp} \nabla)\vct{\mu}}{\rho} - \frac{\adj{\vct{\mu}}^{\tp}\vct{\mu} \vct{\mu}^{\tp}\nabla \rho}{\rho^2} + \adj{\vct{\mu}}^{\tp}\nabla \bar{P} \\
          + \adj{\rho} \divergence\left(\vct{\mu} \right) \\
          + \adj{E} \left(\frac{\vct{\mu}^{\tp} \nabla (E + \bar{P}) + (E + \bar{P}) \divergence(\vct{\mu})}{\rho} - \frac{(E + \bar{P})\vct{\mu}^{\tp} \nabla \rho}{\rho^2}\right)
     \end{bmatrix}.
\end{align}
To evaluate the source term, we differentiate this expression with respect to the undifferentiated occurrences of $\vct{q}$, holding $\Sigma$ and $\cnst{D}\vct{q}$ fixed. 
We denote this operation by $\nabla^{\Sigma, \cnst{D}\vct{q}}_{\vct{q}}$.
By symmetry of the second derivatives, $\divergence\left(\cnst{D}^{\Sigma}_{\vct{q}}\vct{F}\right) = \cnst{D}^{\Sigma, \cnst{D}\vct{q}}_{\vct{q}}\left(\divergence\left(\vct{F}\right)\right)$.
Thus,
\begin{align}
     \adj{\vct{q}}^{\tp}\! \divergence\left(\cnst{D}^{\Sigma}_{\vct{q}}\vct{F}\right)
     \!=\!
     \begin{pmatrix}
          \begin{bmatrix}
          \frac{\divergence(\vct{\mu})\Id   
               +\left(\cnst{D} \vct{\mu} \right)^{\tp}}{\rho}\adj{\vct{\mu}}
          +  \cnst{D} \left(P_{\vct{\mu}}\right) \adj{\vct{\mu}}
          - \frac{\nabla \rho \vct{\mu}^{\tp}
               + \vct{\mu}^{\tp} \nabla \rho}{\rho^2} \adj{\vct{\mu}}\\
          + \left(\frac{\nabla(\bar{H})
               + \cnst{D}\left(P_{\vct{\mu}}\right) \vct{\mu}
               + P_{\vct{\mu}} \divergence(\vct{\mu})}{\rho}
          -  \frac{\bar{H}\nabla \rho
               + P_{\vct{\mu}} \vct{\mu}^{\tp} \nabla \rho}{\rho^2} 
               \right) \adj{E} 
          \end{bmatrix}
          \\
          \begin{bmatrix}
          -\frac{\adj{\vct{\mu}}^{\tp}\divergence(\vct{\mu})\vct{\mu} + \adj{\vct{\mu}}^{\tp}(\vct{\mu}^{\tp} \nabla)\vct{\mu}}{\rho^2} + 2 \frac{\adj{\vct{\mu}}^{\tp}\vct{\mu} \vct{\mu}^{\tp}\nabla \rho}{\rho^3} + \adj{\vct{\mu}}^{\tp}\nabla P_{\rho}\\
          + \adj{E} \left(\frac{\vct{\mu}^{\tp} \nabla P_\rho + P_\rho \divergence(\vct{\mu})}{\rho} 
          - \frac{\vct{\mu}^{\tp} \nabla \bar{H} + \bar{H} \divergence(\vct{\mu})}{\rho^2}
          - \frac{P_\rho\vct{\mu}^{\tp} \nabla \rho}{\rho^2} 
          + 2 \frac{\bar{H}\vct{\mu}^{\tp} \nabla \rho}{\rho^3}
          \right) 
          \end{bmatrix}
          \\
          \adj{\vct{\mu}}^{\tp}\nabla P_E 
          + \adj{E} \left(\frac{\vct{\mu}^{\tp} \nabla P_E + (1 + P_E) \divergence(\vct{\mu})}{\rho} 
               - \frac{(1 + P_E)\vct{\mu}^{\tp} \nabla \rho}{\rho^2}
          \right)
     \end{pmatrix} &\\
     = 
     \begin{pmatrix}
          \begin{bmatrix}
          \left(\frac{\divergence(\vct{\mu})\Id   
               +\left(\cnst{D} \vct{\mu} \right)^{\tp}}{\rho}
          +  \cnst{D} \left(P_{\vct{\mu}}\right)
          - \frac{\nabla \rho \vct{\mu}^{\tp}
               + \vct{\mu}^{\tp} \nabla \rho}{\rho^2} \right) \adj{\vct{\mu}}\\
          + \left(\frac{\nabla(\bar{H})
               + \cnst{D}\left(P_{\vct{\mu}}\right) \vct{\mu}
               + P_{\vct{\mu}}  \divergence(\vct{\mu})}{\rho}
          -  \frac{\bar{H}\nabla \rho
               + \left(\nabla_{\vct{\mu}} P\right) \vct{\mu}^{\tp} \nabla \rho}{\rho^2} 
               \right) \adj{E} 
          \end{bmatrix}
               \\
          \begin{bmatrix}
          \adj{\vct{\mu}}^{\tp} \left(-\frac{\divergence(\vct{\mu})\vct{\mu} + (\vct{\mu}^{\tp} \nabla)\vct{\mu}}{\rho^2} + 2 \frac{\vct{\mu} \vct{\mu}^{\tp}\nabla \rho}{\rho^3} + \nabla P_{\rho} \right)\\
          + \adj{E} \left(\frac{\vct{\mu}^{\tp} \nabla (P_\rho) + P_\rho \divergence(\vct{\mu})}{\rho} 
          - \frac{\vct{\mu}^{\tp} \nabla (\bar{H}) + \bar{H} \divergence(\vct{\mu})}{\rho^2}
          - \frac{P_\rho\vct{\mu}^{\tp} \nabla \rho}{\rho^2} 
          + 2 \frac{\bar{H}\vct{\mu}^{\tp} \nabla \rho}{\rho^3}
          \right) 
          \end{bmatrix}
          \\
          \adj{\vct{\mu}}^{\tp}\nabla P_E 
          + \adj{E} \left(\frac{\vct{\mu}^{\tp} \nabla (P_E) + (1 + P_E) \divergence(\vct{\mu})}{\rho} 
               - \frac{(1 + P_E)\vct{\mu}^{\tp} \nabla \rho}{\rho^2}
          \right)
     \end{pmatrix} &\\
     = 
     \begin{pmatrix}
          \begin{bmatrix}
          \left(\frac{\divergence(\vct{\mu})\Id   
               +\left(\cnst{D} \vct{\mu} \right)^{\tp}}{\rho}
          +  \cnst{D} \left(P_{\vct{\mu}}\right)
          - \frac{\nabla \rho \vct{\mu}^{\tp}
               + \vct{\mu}^{\tp} \nabla \rho}{\rho^2} \right) \adj{\vct{\mu}}\\
          + \left(\frac{\nabla\bar{H}
               + \cnst{D}\left(P_{\vct{\mu}}\right) \vct{\mu}
               + P_{\vct{\mu}} \divergence(\vct{\mu})}{\rho}
          -  \frac{\bar{H}\nabla \rho
               + P_{\vct{\mu}} \vct{\mu}^{\tp} \nabla \rho}{\rho^2} 
               \right) \adj{E} 
          \end{bmatrix}
               \\
          \begin{bmatrix}
          \adj{\vct{\mu}}^{\tp} \left(-\frac{\divergence(\vct{\mu})\vct{\mu} + (\vct{\mu}^{\tp} \nabla)\vct{\mu}}{\rho^2} + 2 \frac{\vct{\mu} \vct{\mu}^{\tp}\nabla \rho}{\rho^3} + \nabla P_{\rho} \right)\\
          + \adj{E} \left(\frac{\vct{\mu}^{\tp} \nabla P_\rho + P_\rho \divergence(\vct{\mu})}{\rho} 
          - \frac{\vct{\mu}^{\tp} \nabla \bar{H} + \bar{H} \divergence(\vct{\mu})}{\rho^2}
          - \frac{P_\rho\vct{\mu}^{\tp} \nabla \rho}{\rho^2} 
          + 2 \frac{\bar{H}\vct{\mu}^{\tp} \nabla \rho}{\rho^3}
          \right) 
          \end{bmatrix}
          \\
          \adj{\vct{\mu}}^{\tp}\nabla P_E
          + \adj{E} \left(\frac{\vct{\mu}^{\tp} \nabla (P_E) + (1 + P_E) \divergence(\vct{\mu})}{\rho} 
               - \frac{(1 + P_E)\vct{\mu}^{\tp} \nabla \rho}{\rho^2}
          \right)
     \end{pmatrix}.& 
\end{align}
Together with the nonlocal term from the advective form, we obtain \cref{eqn:adjoint_abstract_nonlocal_conservative} in the form
{\small
\begin{equation}
\label{eqn:igr_adjoint_conservative}
     \begin{split}
          -
          \partial_t
          \begin{pmatrix}
               \adj{\vct{\mu}}\\
               \adj{\rho}\\
               \adj{E}
          \end{pmatrix}
          -
          \divergence
          \begin{pmatrix}
               \frac{\adj{\vct{\mu}} \vct{\mu}^\tp + \adj{\vct{\mu}}^{\tp} \vct{\mu} \Id}{\rho}
                    + P_{\vct{\mu}} \adj{\vct{\mu}}^{\tp}
                    + \adj{\rho} \Id
                    + \adj{E} \frac{\bar{H}}{\rho} \Id
                    + \adj{E} \frac{P_{\vct{\mu}}\vct{\mu}^{\tp}}{\rho}\\
               - \frac{\adj{\vct{\mu}}^\tp \vct{\mu} \vct{\mu}^\tp}{\rho^2}
                    + P_{\rho} \adj{\vct{\mu}}^{\tp}
                    + \adj{E} \left(\frac{P_{\rho}}{\rho} - \frac{\bar{H}}{\rho^2}\right) \vct{\mu}^{\tp}\\
               P_E \adj{\vct{\mu}}^{\tp} + \adj{E} \frac{\left(1 + P_E \right) \vct{\mu}^{\tp}}{\rho}
          \end{pmatrix}
          -
          \left(\divergence\left(\adj{\vct{\mu}}\right) + \frac{\left[\cnst{D} \adj{E}\right]\vct{\mu}}{\rho}\right)
          \Sigma_{\vct{q}} 
          \\
          +
          \begin{pmatrix}
               \begin{bmatrix}
               \left(\frac{\divergence(\vct{\mu})\Id
                    +\left(\cnst{D} \vct{\mu} \right)^{\tp}}{\rho}
               +  \cnst{D} \left(P_{\vct{\mu}}\right)
               - \frac{\nabla \rho \vct{\mu}^{\tp}
                    + \vct{\mu}^{\tp} \nabla \rho}{\rho^2} \right) \adj{\vct{\mu}}\\
               + \left(\frac{\nabla\bar{H}
                    + \cnst{D}\left(P_{\vct{\mu}}\right) \vct{\mu}
                    + P_{\vct{\mu}} \divergence(\vct{\mu})}{\rho}
               -  \frac{\bar{H}\nabla \rho
                    + P_{\vct{\mu}} \vct{\mu}^{\tp} \nabla \rho}{\rho^2}
                    \right) \adj{E}
               \end{bmatrix}
                    \\
               \begin{bmatrix}
               \adj{\vct{\mu}}^{\tp} \left(-\frac{\divergence(\vct{\mu})\vct{\mu} + (\vct{\mu}^{\tp} \nabla)\vct{\mu}}{\rho^2} + 2 \frac{\vct{\mu} \vct{\mu}^{\tp}\nabla \rho}{\rho^3} + \nabla P_{\rho} \right)\\
               + \adj{E} \left(\frac{\vct{\mu}^{\tp} \nabla P_\rho + P_\rho \divergence(\vct{\mu})}{\rho}
               - \frac{\vct{\mu}^{\tp} \nabla \bar{H} + \bar{H} \divergence(\vct{\mu})}{\rho^2}
               - \frac{P_\rho\vct{\mu}^{\tp} \nabla \rho}{\rho^2}
               + 2 \frac{\bar{H}\vct{\mu}^{\tp} \nabla \rho}{\rho^3}
               \right)
               \end{bmatrix}
               \\
               \adj{\vct{\mu}}^{\tp}\nabla P_E
               + \adj{E} \left(\frac{\vct{\mu}^{\tp} \nabla (P_E) + (1 + P_E) \divergence(\vct{\mu})}{\rho}
                    - \frac{(1 + P_E)\vct{\mu}^{\tp} \nabla \rho}{\rho^2}
               \right)
          \end{pmatrix}
=
          \begin{pmatrix}
               \vct{g}_{\vct{\mu}}\left(\vct{x}, t\right) \\
               g_\rho\left(\vct{x}, t\right) \\
               g_E\left(\vct{x}, t\right)
          \end{pmatrix}
          .
     \end{split}
\end{equation}
}
\subsection{Computation of the nonlocal term}
We close these equations via the $L^2$ gradient $\Sigma_{\vct{q}}$ of $\Sigma$.
The entropic pressure $\Sigma$ is defined by the elliptic equation
\begin{align}
     \frac{\Sigma}{\rho} - \alpha \divergence \left(\frac{\nabla \Sigma}{\rho} \right) 
     &= \alpha \left(\trace^2\left(\cnst{D}\vct{u}\right) + \trace\left(\left(\cnst{D}\vct{u}\right)^2\right)  \right)\\
     &= \alpha \left(\trace^2\left(\cnst{D}\left(\frac{\vct{\mu}}{\rho}\right)\right) + \trace\left(\left(\cnst{D} \left(\frac{\vct{\mu}}{\rho}\right)\right)^2\right)  \right).
\end{align}
We denote the elliptic operator by $\mathcal{L}[\cdot] = \frac{(\cdot)}{\rho} - \alpha \divergence\!\left(\frac{\nabla(\cdot)}{\rho}\right)$.
\begin{align}
     &\int \limits_{\Omega} \left(\divergence\left(\adj{\vct{\mu}}\right) + \frac{\left[\cnst{D} \adj{E}\right]\vct{\mu}}{\rho}\right)
     \Sigma_{\vct{q}} \cdot \sens{\vct{q}} \D \vct{x}\\ 
     =& \left.\frac{\D}{\D \delta} \int \limits_{\Omega} \left(\divergence\left(\adj{\vct{\mu}}\right) + \frac{\left[\cnst{D} \adj{E}\right]\vct{\mu}}{\rho}\right)
     \Sigma\left(\vct{q} + \delta \sens{\vct{q}}\right) \D \vct{x} \right|_{\delta = 0} \\
     =&
     \alpha \int \limits_{\Omega} \left(\divergence\left(\adj{\vct{\mu}}\right) + \frac{\left[\cnst{D} \adj{E}\right]\vct{\mu}}{\rho}\right)
     \mathcal{L}^{-1}
     \Bigg(\trace\left( \cnst{D}\vct{u}\right) \trace\left(2 \cnst{D}\sens{\vct{u}} + \frac{\sens{\rho}}{\rho}\cnst{D} \vct{u}\right) \\ & \qquad \qquad + \trace\left(\cnst{D}\vct{u} \left(2 \cnst{D}\sens{\vct{u}} + \frac{\sens{\rho}}{\rho} \cnst{D}\vct{u}\right)\right) - \nabla\left(\frac{\sens{\rho}}{\rho}\right) \cdot \frac{\nabla \Sigma}{\rho} \Bigg)
     \D \vct{x}\\
     =&
     \alpha \int \limits_{\Omega} \mathcal{L}^{-1} \left(\divergence\left(\adj{\vct{\mu}}\right) + \frac{\left[\cnst{D} \adj{E}\right]\vct{\mu}}{\rho}\right)
     \Bigg(2 \trace\left( \cnst{D}\vct{u}\right) \trace\left( \cnst{D}\sens{\vct{u}}\right) + 2 \trace\left(\cnst{D}\vct{u} \cnst{D}\sens{\vct{u}}\right) \\ 
     & \qquad \qquad \qquad +\left(\trace^2\left( \cnst{D}\vct{u}\right) + \trace(\left[\cnst{D}\vct{u}\right]^2) \right) \frac{\sens{\rho}}{\rho} - \nabla\left(\frac{\sens{\rho}}{\rho}\right) \cdot \frac{\nabla \Sigma}{\rho} \Bigg)
     \D \vct{x},
\end{align}
where we have used the self-adjointness of the solution operator $\mathcal{L}^{-1}$.
Defining $\Pi$ as 
\begin{equation}
     \label{eq:pi_definition}
      \frac{\Pi}{\rho} - \alpha \divergence \left(\frac{\nabla \Pi}{\rho} \right) 
     = \alpha \left(\divergence\left(\adj{\vct{\mu}}\right) + \frac{\left[\cnst{D} \adj{E}\right]\vct{\mu}}{\rho}\right)
\end{equation}
with periodic boundary conditions, we obtain 
\begin{align}
     &\int \limits_{\Omega} \left(\divergence\left(\adj{\vct{\mu}}\right) + \frac{\left[\cnst{D} \adj{E}\right]\vct{\mu}}{\rho}\right)
     \Sigma_{\vct{q}} \cdot \sens{\vct{q}} \D \vct{x} 
     =
     \int \limits_{\Omega} \Pi
     \Bigg(2 \trace\left( \cnst{D}\vct{u}\right) \trace\left( \cnst{D}\sens{\vct{u}}\right) + 2 \trace\left(\cnst{D}\vct{u} \cnst{D}\sens{\vct{u}}\right) \\
     & \qquad \qquad \qquad\qquad \qquad \qquad \qquad +\left(\trace^2\left( \cnst{D}\vct{u}\right) + \trace(\left[\cnst{D}\vct{u}\right]^2) \right) \frac{\sens{\rho}}{\rho} - \nabla\left(\frac{\sens{\rho}}{\rho}\right) \cdot \frac{\nabla \Sigma}{\rho} \Bigg)
     \D \vct{x} \\
     =& 
     \int \limits_{\Omega}  
     -2\left(\nabla\left(\Pi \divergence(\vct{u}) \right) + \divergence\left(\Pi [\cnst{D}\vct{u}]^{\tp}\right) \right) \cdot \sens{\vct{u}} \\
     & \qquad + \left(\frac{\Pi\left(\trace^2\left( \cnst{D}\vct{u}\right) + \trace(\left[\cnst{D}\vct{u}\right]^2) \right)
     + \divergence\left(\Pi \frac{\nabla \Sigma}{\rho}\right)}{\rho}\right) \sens{\rho}
     \D \vct{x} \\
     =& 
     \int \limits_{\Omega}  
     -\frac{2\left(\nabla\left(\Pi \divergence(\vct{u}) \right) + \divergence\left(\Pi [\cnst{D}\vct{u}]^{\tp}\right) \right)}{\rho} \cdot \sens{\vct{\mu}} \\
     &+ \Bigg(\frac{\Pi\left(\trace^2\left( \cnst{D}\vct{u}\right) + \trace(\left[\cnst{D}\vct{u}\right]^2)\right) 
     + \divergence\left(\Pi \frac{\nabla \Sigma}{\rho}\right)}{\rho}\\ 
     & \qquad+ \frac{2\left(\nabla\left(\Pi \divergence(\vct{u}) \right) + \divergence\left(\Pi [\cnst{D}\vct{u}]^{\tp}\right) \right) \cdot \vct{\mu}}{\rho^2}\Bigg) \sens{\rho}
     \D \vct{x}, 
\end{align}
where we have again used the fact that $\sens{\vct{u}} = \frac{\sens{\vct{\mu}}}{\rho} - \frac{\vct{\mu} \sens{\rho}}{\rho^2}$.
Together, this implies 
\begin{equation}
\label{eqn:nonlocal_source_pi}
{\scriptstyle \left(\divergence\left(\adj{\vct{\mu}}\right) + \frac{\left[\cnst{D} \adj{E}\right]\vct{\mu}}{\rho}\right)
\Sigma_{\vct{q}} }
=
\begin{pmatrix}
     - \frac{2\left(\nabla\left(\Pi \divergence(\vct{u}) \right) + \divergence\left(\Pi [\cnst{D}\vct{u}]^{\tp}\right) \right)}{\rho} \\
     \frac{\Pi\left(\trace^2\left( \cnst{D}\vct{u}\right) + \trace(\left[\cnst{D}\vct{u}\right]^2) \right)
     + \divergence\left(\Pi \frac{\nabla \Sigma}{\rho}\right)}{\rho} 
     + \frac{2\left(\nabla\left(\Pi \divergence(\vct{u}) \right) + \divergence\left(\Pi [\cnst{D}\vct{u}]^{\tp}\right) \right) \cdot \vct{\mu}}{\rho^2}\\
     0
\end{pmatrix}.
\end{equation}
\section{Numerical Implementation}
The equations of IGR and its sensitivities are not restricted to a particular discretization scheme. 
For the purpose of this work, we use a discontinuous Galerkin (DG) method for their spatial discretization. 
This automatically provides the necessary stabilization required to solve the transport-dominant forward and adjoint problems as well as grid function access to the entropic pressure. 
The latter aids in evaluating the $\nabla \Sigma$ appearing in the adjoint PDE.
\subsection{Discretization of the forward problem}
\label{sec:forward_discretization}
We discretize IGR and its forward sensitivities by a nodal discontinuous Galerkin (DG) method \cite{cockburn1998local,cockburn2012discontinuous} introduced for IGR by \cite{eyob2026discontinuous}, to which we refer for a detailed description.
On a conforming partition $\mathcal{T}_h = \{K\}$ of $\Omega$ with faces $\mathcal{F}_h$, all unknowns are approximated in the space $V_h^p$ of piecewise polynomials of degree at most $p$. Interpolation and quadrature uses tensor-product Gauss--Lobatto points.
We use explicit strong-stability-preserving Runge--Kutta methods for time integration, with an elliptic solve at every stage \cite{shu1988total,hadjimichael2013strong}.

\subsubsection{DG discretization of the hyperbolic part}
The hyperbolic part uses the standard weak form DG with the local Lax--Friedrichs/Rusanov numerical flux
\begin{equation}
\label{eqn:numerical_flux}
     \widehat{\vct{f}}_{\vct{n}}(\vct{q}^-, \vct{q}^+, \Sigma^-, \Sigma^+)
     = \frac{1}{2}\left(\vct{F}(\vct{q}^-, \Sigma^-)\vct{n} + \vct{F}(\vct{q}^+, \Sigma^+)\vct{n}\right)
     - \frac{\lambda_f}{2}(\vct{q}^+ - \vct{q}^-),
\end{equation}
where $\lambda_f = \max(|\vct{u}^- \cdot \vct{n}| + c^-,\, |\vct{u}^+ \cdot \vct{n}| + c^+)$ is an upper bound on the maximum wave speed at the face, with $c^{\pm}$ denoting the sound speeds.
The same formulation applies to the forward sensitivity equations, with additional flux dependence on $\sens{\Sigma}$.

\subsubsection{SIP discretization of the elliptic equation}
The entropic pressure $\Sigma$ and its sensitivity $\sens{\Sigma}$ solve elliptic PDEs of the common form $\Sigma/\rho - \alpha \divergence\left(\rho^{-1}\nabla \Sigma\right) = R(\vct{q})$, which we discretize with the symmetric interior penalty (SIP) method \cite{douglas2008interior,wheeler1978elliptic,arnold1982interior}: find $\Sigma \in V_h^p$ such that $a_h(\Sigma, \phi) = \ell_h(\phi)$ for all $\phi \in V_h^p$, with
\begin{align}
\label{eqn:dg_elliptic}
     a_h(\Sigma, \phi) &=
     \sum_{K \in \mathcal{T}_h} \int_K \frac{\Sigma \phi}{\rho} \D \vct{x}
     + \alpha \sum_{K \in \mathcal{T}_h} \int_K \frac{\nabla \Sigma \cdot \nabla \phi}{\rho} \D \vct{x} \notag \\
     &\quad - \alpha \sum_{f \in \mathcal{F}_h} \int_f \avg{\frac{\nabla \Sigma}{\rho}} \cdot \vct{n}_f \jump{\phi} \D s
     - \alpha \sum_{f \in \mathcal{F}_h} \int_f \avg{\frac{\nabla \phi}{\rho}} \cdot \vct{n}_f \jump{\Sigma} \D s \notag \\
     &\quad + \alpha \sum_{f \in \mathcal{F}_h} \int_f \frac{\eta_p}{h_f} \avg{\frac{1}{\rho}} \jump{\Sigma} \jump{\phi} \D s, \qquad
     \ell_h(\phi) = \sum_{K \in \mathcal{T}_h} \int_K R(\vct{q}) \phi \D \vct{x}.
\end{align}
Here, $\jump{\cdot}$ and $\avg{\cdot}$ denote the jump and average of the traces across a face $f$ with unit normal $\vct{n}_f$ and diameter $h_f$, and $\eta_p > 0$ is a penalty parameter chosen large enough to ensure coercivity.
The resulting linear systems are symmetric positive definite and can be solved by iterative methods, such as diagonally preconditioned conjugate gradient.

\begin{remark}[Elliptic solvers for differentiation]
     We found that automatic differentiation needs more accurate elliptic solves in the forward problem and that differentiation through preconditioned conjugate gradient (PCG) is unstable. 
     Thus, we use a large number of Jacobi iterations when using forward-mode automatic differentiation. 
     In reverse mode automatic differentiation, this costs too much memory.
     Thus, we use a smaller number of Chebyshev iterations for reverse-mode automatic differentiation.
     This could be avoided by implementing custom sensitivities for the elliptic solve.
\end{remark}

\subsection{Discretization of the adjoint problem}
By storing the states of the forward evolution as grid functions, we obtain the coefficients for the adjoint problem, which we integrate backwards in time. 
Checkpointing reduces the memory requirements by storing only a subset of the forward states and recomputing the rest as needed \cite{charpentier2001checkpointing,griewank2008evaluating,wang2009minimal}.
To simplify the implementation, we focus on problem sizes that do not require checkpointing on a laptop. 

\subsubsection{Conservative form and time reversal}
The conservative form \cref{eqn:adjoint_abstract_nonlocal_conservative} of the adjoint has a local flux and a nonlocal source term.
Reversing time by setting $\tau = T - t$, the adjoint equation becomes a forward-in-$\tau$ conservation law with source,
\begin{equation}
\label{eqn:adjoint_conservative}
     \partial_\tau \adj{\vct{q}} + \divergence(\adj{\vct{F}}_{\mathrm{loc}}) = \vct{g} + \vct{S}_{\mathrm{loc}} + \vct{S}_{\Sigma},
\end{equation}
where $\adj{\vct{F}}_{\mathrm{loc}}$ is the local adjoint flux (depending only on $P$ derivatives), $\vct{S}_{\mathrm{loc}}$ is the local source term, and $\vct{S}_{\Sigma}$ is the nonlocal source arising from the entropic pressure.

\subsubsection{Derivatives of the physical pressure}
The local adjoint flux and source require derivatives of the physical pressure $P$.
The equation of state depends on the density and specific internal energy $e = E/\rho - |\vct{u}|^2/2$ and the chain rule yields
\begin{equation}
     \partial_E P = \frac{P_e}{\rho}, \qquad
     \nabla_{\vct{\mu}} P = -\frac{P_e \vct{u}}{\rho}, \qquad
     \partial_\rho P = P_\rho + P_e\!\left(-\frac{E}{\rho^2} + \frac{|\vct{\mu}|^2}{\rho^3}\right),
\end{equation}
where $P_\rho = \partial P/\partial \rho$ and $P_e = \partial P/\partial e$ are partial derivatives of the equation of state.
For a polytropic gas $P = (\gamma - 1)\rho e$, we have $P_\rho = (\gamma-1)e$ and $P_e = (\gamma-1)\rho$, yielding
\begin{equation}
\label{eqn:polytropic_pressure_derivs}
     \partial_E P = \gamma - 1, \qquad
     \nabla_{\vct{\mu}} P = -(\gamma-1)\vct{u}, \qquad
     \partial_\rho P = \frac{(\gamma-1)|\vct{u}|^2}{2}.
\end{equation}

\subsubsection{Spatial derivatives for the local source term}
The local adjoint flux $\adj{\vct{F}}_{\mathrm{loc}}$ is a pointwise function of $\adj{\vct{q}}$, $\vct{q}$, and the pressure derivatives $P_E$, $P_\rho$, $P_{\vct{\mu}}$.
The local source term $\vct{S}_{\mathrm{loc}}$ needs spatial derivatives of state variables and pressure (derivatives).
The derivatives $\cnst{D}\vct{\mu}$, $\divergence(\vct{\mu}) = \trace(\cnst{D}\vct{\mu})$, $\nabla\rho$, and $\nabla E$ are computed element-locally by differentiating grid functions.
The spatial gradient $\nabla\bar{H}$ (with $\bar{H} = E + P + \Sigma$) uses $\nabla\bar{H} = \nabla E + \nabla P + \nabla\Sigma$, for $\nabla\Sigma$ computed element-locally from the $\Sigma$ grid function.
The remaining spatial derivatives needed are $\nabla P$, $\nabla(P_\rho)$, $\nabla(P_E)$, and $\cnst{D}(P_{\vct{\mu}})$.
For the polytropic equation of state $P = (\gamma - 1)(E - |\vct{\mu}|^2/(2\rho))$, these are
\begin{align}
\label{eqn:polytropic_spatial_derivs}
     \nabla P &= (\gamma - 1)\!\left(\nabla E - \frac{(\cnst{D}\vct{\mu})^{\tp}\vct{\mu}}{\rho} + \frac{|\vct{\mu}|^2}{2\rho^2}\nabla\rho\right), \\
     \nabla(P_\rho) &= (\gamma - 1)\!\left(\frac{(\cnst{D}\vct{\mu})^{\tp}\vct{\mu}}{\rho^2} - \frac{|\vct{\mu}|^2}{\rho^3}\nabla\rho\right), \\
     \nabla(P_E) &= 0, \\
     \cnst{D}(P_{\vct{\mu}}) &= -(\gamma - 1)\cnst{D}\vct{u} = -(\gamma - 1)\!\left(\frac{\cnst{D}\vct{\mu}}{\rho} - \frac{\vct{\mu}(\nabla\rho)^{\tp}}{\rho^2}\right).
\end{align}
All right-hand side quantities are element-local derivatives of state variables $\vct{\mu}$, $\rho$, $E$.

\subsubsection{Computation of the nonlocal source term}
The nonlocal source term $\vct{S}_{\Sigma}$ involves the $L^2$ gradient $\Sigma_{\vct{q}}$.
The definition \cref{eqn:nonlocal_source_pi} of $\Sigma_{\vct{q}}$ in terms of the $\Pi$ defined in \cref{eq:pi_definition} only needs a single elliptic solve.
Discretely, we define $\Pi \in V_h^p$ by
\begin{equation}
\label{eqn:dg_adjoint_elliptic}
     a_h(\Pi, \psi) = \alpha \sum_{K \in \mathcal{T}_h} \int_K \left(\divergence\!\left(\adj{\vct{\mu}}\right) + \frac{\left[\cnst{D}\adj{E}\right]\vct{\mu}}{\rho}\right) \psi \D \vct{x} \quad \text{for all } \psi \in V_h^p,
\end{equation}
where $a_h$ is the same SIP bilinear form from \cref{eqn:dg_elliptic}.
Since $a_h$ is symmetric, $\Pi$ can be computed using the same stiffness matrix as the forward $\Sigma$ solve.
The spatial derivatives of $\Pi$ (appearing through $\nabla(\Pi\divergence(\vct{u}))$ and $\divergence(\Pi[\cnst{D}\vct{u}]^{\tp})$) are computed element-locally by differentiating the polynomial expansion of $\Pi$.

\subsubsection{Semi-discrete DG formulation of the adjoint equation}
\label{sec:dg_adjoint_formulation}
The pressure derivatives, forward solution $\vct{q}$, its entropic pressure $\Sigma$, and the adjoint elliptic field $\Pi$ allow computing the local adjoint flux $\adj{\vct{F}}_{\mathrm{loc}}$, local source $\vct{S}_{\mathrm{loc}}$, and nonlocal source $\vct{S}_{\Sigma}$ element-locally.
The semi-discrete DG formulation reads: find $\adj{\vct{q}}(\tau) \in [V_h^p]^{d+2}$ such that for all $K \in \mathcal{T}_h$ and all $\vct{\phi} \in [V_h^p]^{d+2}$,
\begin{equation}
\label{eqn:dg_adjoint}
\begin{split}
     \int_K \partial_\tau \adj{\vct{q}} \cdot \vct{\phi} \D \vct{x}
     - \int_K \adj{\vct{F}}_{\mathrm{loc}} \mathrel{\colon} \cnst{D}\vct{\phi} \D \vct{x}
     &+ \int_{\partial K} \widehat{\adj{\vct{f}}}_{\vct{n}}(\adj{\vct{q}}^-, \adj{\vct{q}}^+, \vct{q}^{\pm}, \Sigma^{\pm}) \cdot \vct{\phi} \D s\\
     & \qquad \qquad \qquad = \int_K (\vct{g} + \vct{S}_{\mathrm{loc}} + \vct{S}_{\Sigma}) \cdot \vct{\phi} \D \vct{x},
\end{split}
\end{equation}
where $\adj{\vct{F}}_{\mathrm{loc}} = \adj{\vct{F}}_{\mathrm{loc}}(\adj{\vct{q}}, \vct{q}, \Sigma)$, $\vct{S}_{\mathrm{loc}} = \vct{S}_{\mathrm{loc}}(\adj{\vct{q}}, \vct{q}, \Sigma)$, and $\vct{S}_{\Sigma} = \vct{S}_{\Sigma}(\adj{\vct{q}}, \vct{q}, \Sigma, \Pi)$.
The numerical flux mirrors the local Lax--Friedrichs flux \cref{eqn:numerical_flux} of the forward problem,
\begin{equation}
\label{eqn:adjoint_numerical_flux}
     \widehat{\adj{\vct{f}}}_{\vct{n}}(\adj{\vct{q}}^-, \adj{\vct{q}}^+, \vct{q}^{\pm}, \Sigma^{\pm})
     = \frac{1}{2}\left(\adj{\vct{F}}_{\mathrm{loc}}^- \vct{n} + \adj{\vct{F}}_{\mathrm{loc}}^+ \vct{n}\right) - \frac{\lambda_f}{2}(\adj{\vct{q}}^+ - \adj{\vct{q}}^-),
\end{equation}
for $\adj{\vct{F}}_{\mathrm{loc}}^{\pm} = \adj{\vct{F}}_{\mathrm{loc}}(\adj{\vct{q}}^{\pm}, \vct{q}^{\pm}, \Sigma^{\pm})$ and $\lambda_f$ the same wave speed estimate as in \cref{eqn:numerical_flux}.

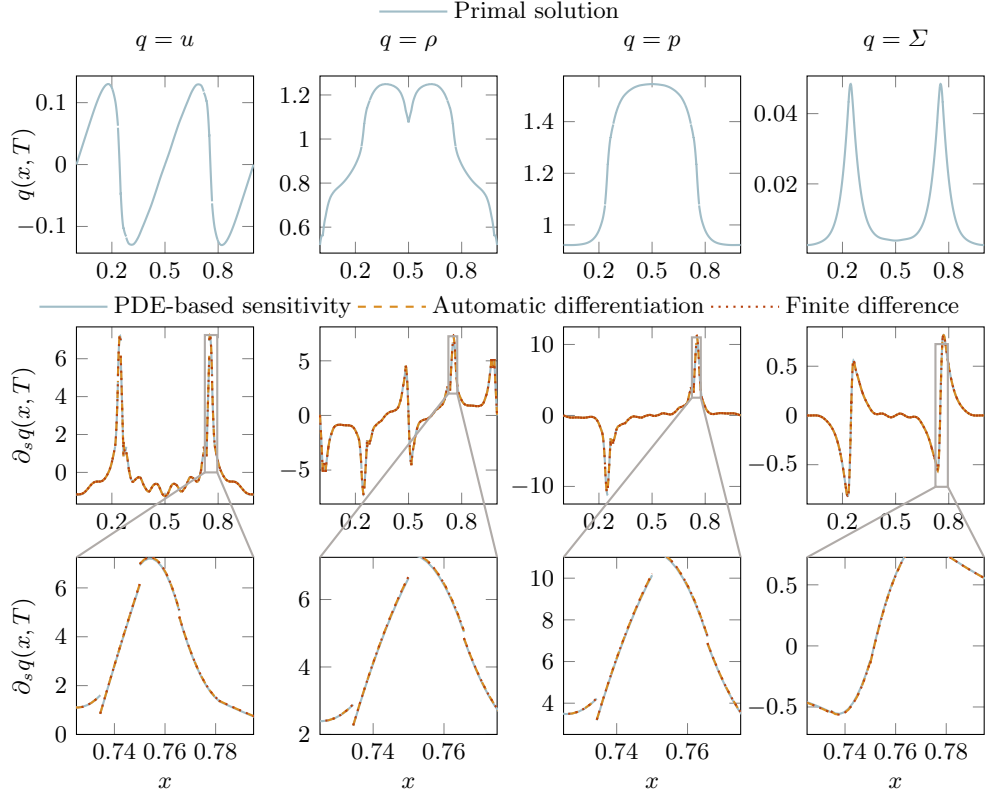
\begin{figure}
     \centering
     \tikzsetnextfilename{shift-sensitivity-postshock}
     \begin{tikzpicture}
     \end{tikzpicture}
     \vspace{-6mm}
     \caption{Velocity sine wave ($A = 1.5$) at $T = 1.7$ (post-shock).
     Top: primal solution showing developed shocks.
     Bottom: sensitivities computed by four methods.
     AD and FD agree (both linearize the discrete operator). 
     The PDE sensitivity differs only slightly, near the regularized shocks.}
     \label{fig:shift_postshock}
\end{figure}

\section{Numerical Experiments}
The code reproducing the following results is available at \url{https://github.com/f-t-s/IGR_sensitivities}.
\subsection{Forward sensitivities}
We begin by demonstrating the computation of forward sensitivities with respect to cyclical shifts and the adiabatic constant $\gamma$.
We compare PDE-based sensitivities via \cref{sec:forward_sensitivities}, automatic differentiation via ForwardDiff \cite{revels2016forward} (AD), and finite differences (FD).

\paragraph{Derivative under cyclic shift}
We consider the initial condition
\begin{equation}
     \label{eqn:sine_wave}
     \rho_0 = 1, \quad u_0(x) = A \sin\!\Big(\frac{2\pi(x - s)}{L}\Big), \quad p_0 = 1,
\end{equation}
with amplitude $A = 1.5$ and period $L = 1$, and compute the sensitivity of the solution with respect to the shift parameter $s$ at $s = 0$.
We use polynomial order $p = 2$, $N_e = 64$ elements, IGR parameter $\alpha = (3 L/N_e)^2$, and FD step size $\varepsilon = 10^{-5}$.
\Cref{fig:shift_postshock} shows results at $T = 1.7 \gg t_s$, well after shock formation.
AD and FD agree with each other and show a minimal discretize-then-differentiate gap to the PDE-based sensitivity.
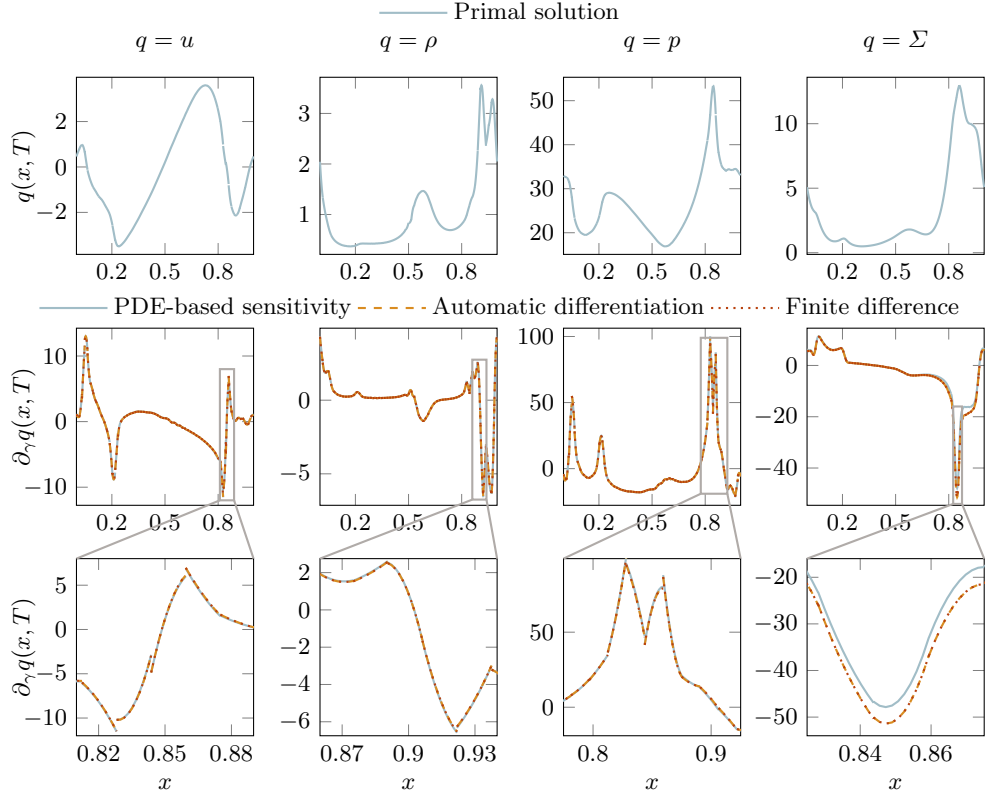
\begin{figure}
     \centering
     \tikzsetnextfilename{blast-gamma-sensitivity}
     \begin{tikzpicture}
          
     \end{tikzpicture}
     \vspace{-8mm}
     \caption{\textbf{Interacting Sedov-like blasts.} Sensitivity $\partial_\gamma q$ at $T = 0.05$.
     Top: primal solution of interacting blast waves.
     Middle: sensitivities $\partial_\gamma u$, $\partial_\gamma \rho$, $\partial_\gamma p$.
     Bottom: blast interaction region.}
     \label{fig:blast_gamma}
\end{figure}

\paragraph{Derivative with respect to the adiabatic constant}
We consider the sensitivity with respect to the adiabatic constant $\gamma$ at $\gamma = 1.4$ of two interacting blast waves using 
\begin{equation}
     \rho_0 = 1, \quad u_0 = 0, \quad p_0(x) = p_{\mathrm{bg}} + p_1 \exp\!\left(-\frac{(x - x_1)^2}{\sigma^2}\right) + p_2 \exp\!\left(-\frac{(x - x_2)^2}{\sigma^2}\right),
\end{equation}
with $p_{\mathrm{bg}} = 1$, $p_1 = 100$, $p_2 = 50$, $x_1 = 0.3$, $x_2 = 0.7$, and $\sigma = 0.1$.
We use polynomial order $p = 2$, $N_e = 64$ elements, IGR parameter $\alpha = 10 (L/N_e)^2$, and FD step size $\varepsilon = 10^{-6}$.
%
\Cref{fig:blast_gamma} shows results at $T = 0.05$, after the two blast waves have interacted.
The PDE-based sensitivity, AD, and FD agree in smooth regions and show a small discretize-then-differentiate gap near shocks.

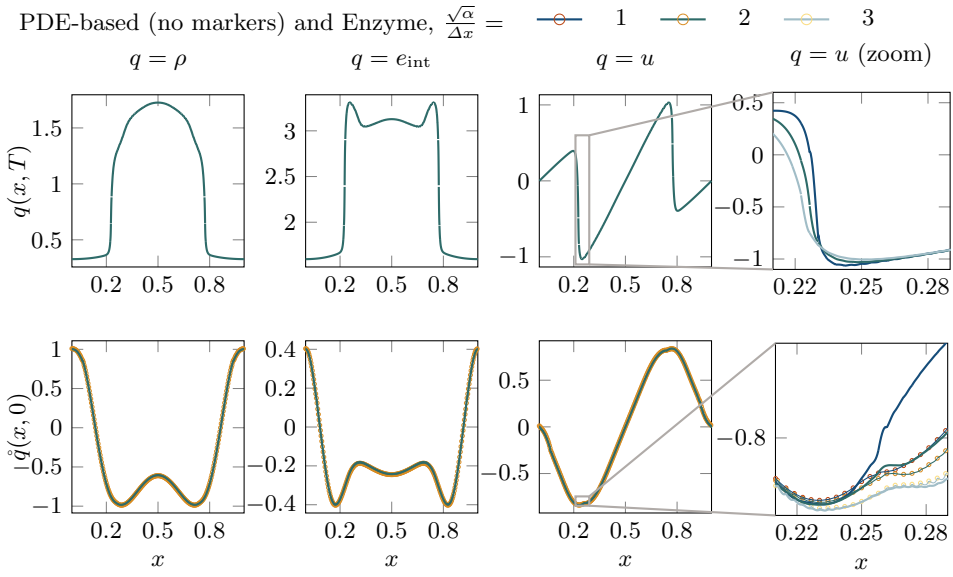
\begin{figure}[htbp]
     \centering
     \tikzsetnextfilename{one-d-adjoints}
     \begin{tikzpicture}
          
     \end{tikzpicture}
     \vspace{-8mm}
     \caption{\textbf{Discrete and continuous adjoints.} We observe strong agreement between adjoints computed by solving \cref{eqn:igr_adjoint_conservative} and by automatic differentiation through the discretization of \cref{eqn:igr_euler}.}
     \label{fig:pde_vs_enzyme_adjoint}
\end{figure}

\begin{figure}
     \centering
     \tikzsetnextfilename{pde-mesh-sweep-merged}
     \begin{tikzpicture}
          
     \end{tikzpicture}
     \vspace{-2mm}
     \caption{\textbf{Scaling across resolutions.} Adjoint solutions across mesh resolutions in two scaling regimes. 
     Top: \emph{Linear} scaling
     $\sqrt{\alpha} = s\, \Delta x$, which keeps the number of cells per layer width $\sqrt{\alpha}$ fixed at $s$.
     Bottom: \emph{Sublinear} scaling $\sqrt{\alpha} = s\, \Delta x_0^{1/3} \Delta x^{2/3}$
     (with $\Delta x_0 = L/128$ to ensure identical $\alpha$ at the coarsest resolution), under which the number of cells per layer grows like $s\,(N_e/128)^{1/3}$.}
     \label{fig:pde_mesh_sweep}
\end{figure}
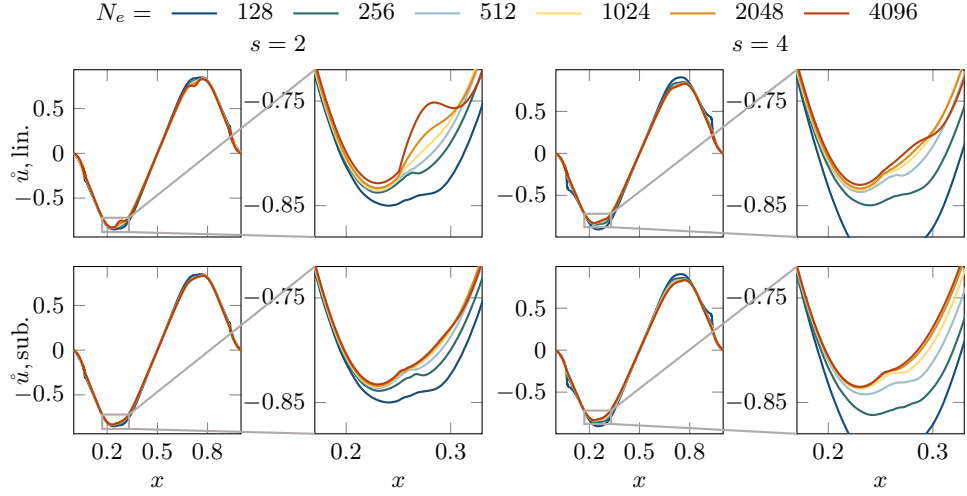

\subsection{Reverse-mode/Adjoint sensitivities}
We now compute numerical solutions of the continuous adjoint equations \cref{eqn:igr_adjoint_conservative} and compare them to automatic differentiation of the forward solve via Enzyme.jl \cite{enzyme1,enzyme2} (for unidimensional problems) or finite difference approximations (for two-dimensional problems).

\paragraph{Discrete and continuous adjoints} As a first step, we consider the initial condition  \cref{eqn:sine_wave}, end time $T = 0.25$, and the terminal linear functional
\begin{equation}
     \label{eqn:weighted_momentum}
     J = \int_\Omega \mu(x, T) \sin\left(2\pi x/L\right) \D x.
\end{equation}
We use $N_e = 256$ elements of polynomial order $p = 3$, and different ratios $\sqrt{\alpha} / \Delta x \in \{1, 2, 3\}$.
We compare the discrete adjoint computed by automatic differentiation with Enzyme to the DG solution of \cref{eqn:igr_adjoint_conservative}, as described in \cref{sec:dg_adjoint_formulation}.
\cref{fig:pde_vs_enzyme_adjoint} shows that discrete and continuous adjoint are in good agreement, which improves as the $\sqrt{\alpha} / \Delta x$ increases. 
This is expected, since this ratio measures the number of elements used to resolve the regularized shock.

\begin{figure}[htbp]
     \centering
     \tikzsetnextfilename{adjoint-convergence-ablation}
     \begin{tikzpicture}
          
     \end{tikzpicture}
     \vspace{-6mm}
     \caption{\textbf{Linear versus sublinear scaling.}
     $L^1$ difference between solutions at consecutive resolutions of the mesh sweep,
     for the primal velocity $u(\cdot, T)$ (left) and the adjoint velocity $\adj{u}(\cdot, 0)$ (right).
     Dashed: linear scaling $\sqrt{\alpha} = s\,\Delta x$; solid: sublinear scaling
     $\sqrt{\alpha} = s\,\Delta x_0^{1/3}\Delta x^{2/3}$; each pair shares its $\alpha$ at $N_e = 128$.
     The primal contracts along both paths; the adjoints of the linear designs break upward at the finest
     resolutions, while their sublinear partners retain convergence.}
     \label{fig:adjoint_convergence_ablation}
\end{figure}
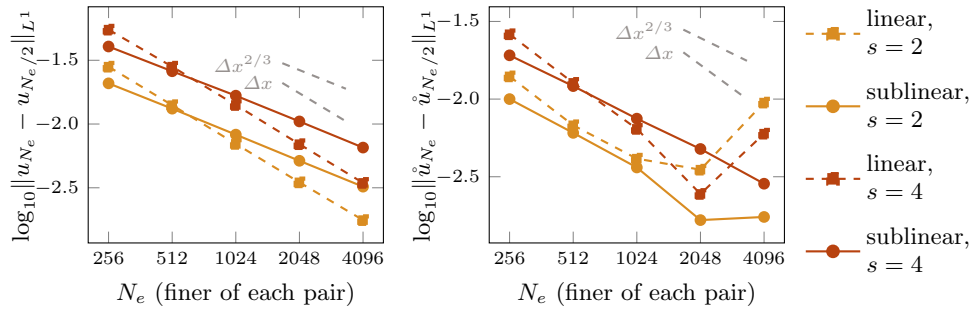

\paragraph{Scaling the regularization with the mesh}
The practical use of IGR requires jointly letting $\alpha$ and $\Delta x$ tend to zero. 
Most applications use $\sqrt{\alpha} \propto \Delta x$ to ensure that the width of the IGR-regularized shock profile is proportional to the grid spacing.
We probe this choice with a mesh sweep over $N_e \in \{128, 256, \ldots, 4096\}$ at polynomial order $p = 4$ and end time $T = 0.25$, starting from the initial condition \cref{eqn:sine_wave} and using the functional \cref{eqn:weighted_momentum}. 
\cref{fig:pde_mesh_sweep,fig:adjoint_convergence_ablation} show that under this scaling, the discretized adjoint solution appears to converge at first, but then continues to change.
A similar effect was studied by \cite{giles2010convergence1,giles2010convergence2} in the context of viscous regularizations.
They show in their setting that scaling the viscosity (and thus the viscous shock width) as $\Delta x^{2/3}$ ensures convergence.
Motivated by their findings, we observe in \cref{fig:pde_mesh_sweep,fig:adjoint_convergence_ablation} that $\sqrt{\alpha} \propto \Delta x^{2/3}$ appears to yield discrete convergence of the continuous IGR adjoint. 
A rigorous analysis of this phenomenon is the subject of future work. 

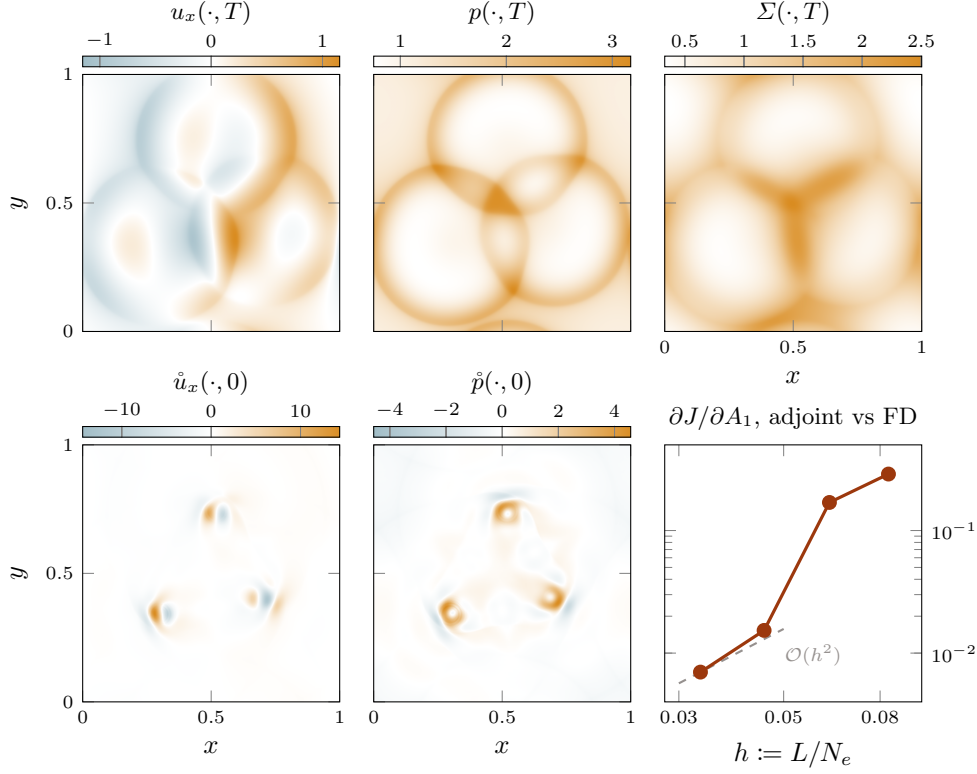
\begin{figure}[htbp]
     \centering
     \tikzsetnextfilename{sedov-heatmaps-2d}
     \begin{tikzpicture}
          
     \end{tikzpicture}
     \vspace{-6mm}
     \caption{\textbf{2D triple Sedov blast} at $N_e = 160$, $p = 4$,
     $T = 0.1$.
     Left: velocity $u_x(\cdot, T)$ and its continuous PDE
     adjoint $\adj{u}_x(\cdot, 0)$ for the mass-averaged kinetic-energy
     objective; center: pressure $p(\cdot, T)$ and adjoint
     $\adj{p}(\cdot, 0)$; right: entropic pressure $\Sigma(\cdot, T)$
     and the mismatch of PDE-adjoint and FD for the directional derivative
     $\D J / \D A_1$ under mesh refinement at fixed $\alpha$.}
     \label{fig:sedov_heatmaps_2d}
\end{figure}

\paragraph{Triple Sedov blast}
We next compute the 2d adjoints on a triple Sedov blast, the superposition of an ambient gas at rest with three Gaussian pressure distributions,
\begin{equation}
     \rho_0 = 1, \quad \vct{u}_0 = \zerovct, \quad
     p_0(\vct{x}) = p_{\mathrm{bg}} + \sum_{k=1}^{3} A_k
       \exp\!\left(-\frac{|\vct{x} - \vct{x}_k|^2}{\sigma^2}\right),
\end{equation}
with $p_{\mathrm{bg}} = 1$, amplitudes $(A_1, A_2, A_3) = (6, 5, 5.5)$,
centers 
$\vct{x}_k = \left(\begin{smallmatrix} 0.30\\ 0.34\end{smallmatrix}\right), 
               \left(\begin{smallmatrix} 0.70\\ 0.40\end{smallmatrix}\right), 
               \left(\begin{smallmatrix} 0.52\\ 0.74\end{smallmatrix}\right)$, and width $\sigma = 0.09$.
The first five panels of \cref{fig:sedov_heatmaps_2d} show the forward solution at $T = 0.1$ and the PDE-based adjoint with respect to the initial condition, for the functional $J$ given by the mass-averaged kinetic energy at $T$.
Both two-dimensional experiments use $N_e = 160$ elements of order $p = 4$ and the fixed regularization $\alpha = 3 (L/20)^2$, which is held constant for the mesh refinement studies, so that refinement targets a single continuous problem.
To verify the PDE-based adjoint, we use it to compute the directional derivative $\D J / \D A_1$ as $\langle \adj{E}(\cdot, 0),\, \partial E_0 / \partial A_1 \rangle$.
We compare this value against a central finite difference approximation and report the result in the last panel of \cref{fig:sedov_heatmaps_2d}.
This discretize-then-differentiate gap decreases with the mesh size.

\begin{figure}[htbp]
     \centering
     \tikzsetnextfilename{oseen-heatmaps-2d}
     \begin{tikzpicture}
          
     \end{tikzpicture}
     \vspace{-6mm}
     \caption{\textbf{Blast--vortex interaction} at $N_e = 160$,
     $p = 4$, $T = 0.1$.
     Left: velocity $u_x(\cdot, T)$ and its continuous PDE adjoint
     $\adj{u}_x(\cdot, 0)$ for the windowed kinetic-energy objective; 
     center: velocity $u_y(\cdot, T)$ and adjoint $\adj{u}_y(\cdot, 0)$;
     right: entropic pressure $\Sigma(\cdot, T)$ and the mismatch of PDE-adjoint and FD for the directional derivative
     $\D J / \D A$ under mesh refinement at fixed $\alpha$.
     The dashed circle on the primal velocity panels marks the $1/e$
     contour of the terminal objective's Gaussian window.}
     \label{fig:oseen_heatmaps_2d}
\end{figure}
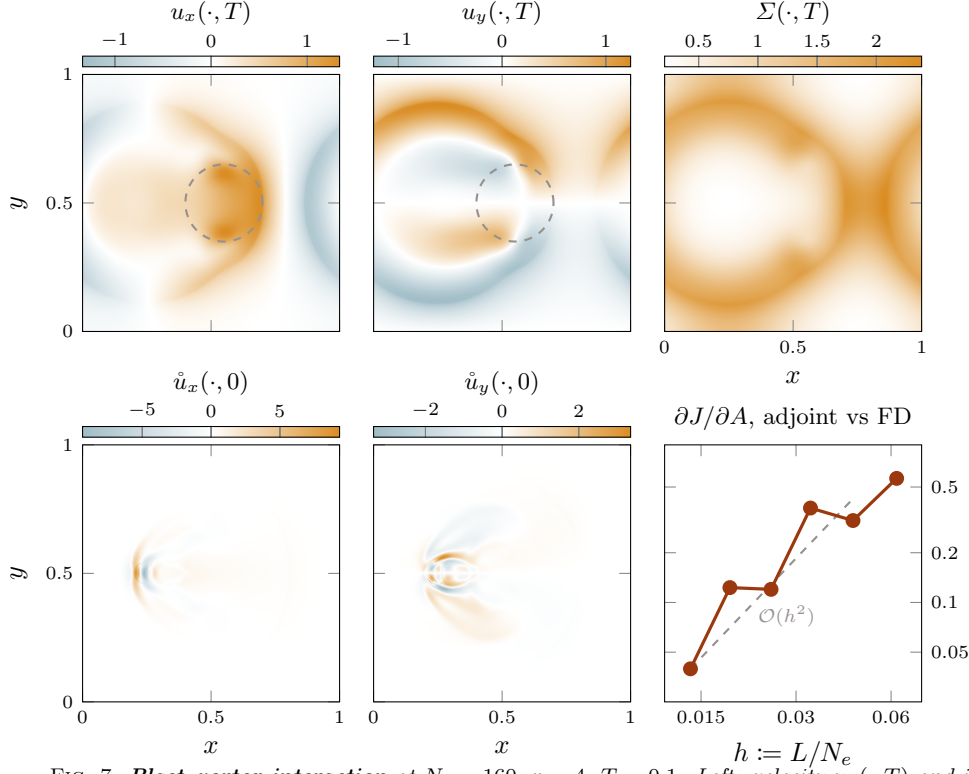%
\paragraph{Blast--vortex interaction}
We now consider the interaction of a Gaussian blast of amplitude $A = 12$ with a pair of counter-rotating Oseen vortices \cite{colonius1991free}.
Each vortex carries Gaussian vorticity with tangential velocity
$v_\theta(r) = \Gamma/(2\pi r)\,\bigl(1 - e^{-r^2/a^2}\bigr)$, core radius
$a = 0.05$, and vortex Mach number $v_{\max}/a_\infty = 0.5$. 
Density and pressure are chosen to ensure radial equilibrium.
The objective $J$ is the kinetic energy in a Gaussian window at the final time $T=0.1$.
As in the previous example, we verify the PDE-based adjoint by computing the resulting directional derivative with respect to $A$ and comparing it to a finite difference approximation. 
\cref{fig:oseen_heatmaps_2d} shows the results.

\section*{Acknowledgments}
This work was supported by the Air Force Office of Scientific Research under award number FA9550-23-1-0668 (Information Geometric Regularization for Simulation and Optimization of Supersonic Flow),
The Predictive Science Academic Alliance Program (PSAAP Award DE-NA0004261 - “The Center for Information Geometric Mechanics and Optimization
(CIGMO)”) managed by the NNSA (National
Nuclear Security Administration) Office of Advanced Simulation, and the Alfred P. Sloan Foundation via a Sloan Research Fellowship in Mathematics.
I thank Matthew Sumanen for his contributions to preliminary work that derived the unidimensional forward sensitivity equation of IGR and Spencer Bryngelson for first making me aware of the challenges due to adjoint computations in flows with shocks. 
I used Claude Code for the following tasks: Checking mathematical derivations, numerical implementation, drafting a first version of section 4, running the numerical experiments, creating TikZ figures, and proofreading of the manuscript. 
\bibliography{references}
\bibliographystyle{plain}

\appendix

\end{document}